\input amstex
\documentstyle{amsppt}
\NoBlackBoxes \magnification=\magstep1 \vsize=22 true cm \hsize=16.1
true cm \voffset=1 true cm
%%%%%%%%%%%%%%%%%%%%%%%%%%%%%%%%%%%%%%%%%%%%%%%%%%%%%%%%%%%%

\define\px{{\partial}_z}
\define\py{{\partial}_{\bar{z}}}
\redefine\id{\operatorname{id}}
\redefine\BMO{\operatorname{BMO}}
\redefine\VMO{\operatorname{VMO}}
\redefine\WP{\operatorname{WP}}
\redefine\QS{\operatorname{QS}}

\redefine\Diff{\operatorname{Diff}}
\redefine\Rot{\operatorname{Rot}}
\topmatter
\title     Weil-Petersson Teichm\"uller space
\endtitle
\title     Weil-Petersson Teichm\"uller space
\endtitle
\author   Shen Yuliang
\endauthor
\affil (Department of Mathematics, Soochow  University)
\endaffil
\address Department of Mathematics,  Soochow University,   Suzhou
215006,  P.\,R. China
\endaddress
\email\nofrills\eighttt  ylshen\@suda.edu.cn
\endemail
\abstract The paper presents some recent results on the Weil-Petersson geometry  theory of the universal Teichm\"uller space, a topic which  is important in Teichm\"uller theory and has wide applications to various areas such as mathematical physics, differential equation and computer vision.

\noindent (1) It is shown that a sense-preserving  homeomorphism $h$ on the unit circle belongs to the Weil-Petersson class, namely, $h$ can be extended to a quasiconformal mapping to the unit disk whose Beltrami coefficient is square integrable in the Poincar\'e metric if and only if $h$ is absolutely continuous  such that $\log h'$ belongs to the Sobolev  class $H^{\frac 12}$. This solves an open problem posed by Takhtajan-Teo [TT2] in 2006 and investigated later by Figalli [Fi],  Gay-Balmaz-Marsden-Ratiu ([GMR], [GR])  and others.

\noindent The intrinsic characterization (1) of the Weil-Petersson class has the following applications which are also explored in this paper:

\noindent (2) It is proved  that there exists a  quasisymmetric homeomorphism  of the Weil-Petersson class which belongs neither to the Sobolev class $H^{\frac 32}$ nor to the Lipschitz class $\Lambda^1$, which was conjectured very recently by Gay-Balmaz-Ratiu [GR] when studying the classical Euler-Poincar\'e equation in the new setting that the involved sense-preserving  homeomorphisms on the unit circle belong to the Weil-Petersson class.

\noindent (3) It is proved that the flows of the $H^{\frac 32}$ vector fields on the unit circle  are contained in the Weil-Petersson class, which was also conjectured by Gay-Balmaz-Ratiu [GR] during their above mentioned research.

\noindent (4) A new  metric is introduced  on the Weil-Petersson Teichm\"uller space and is shown to be topologically equivalent to the Weil-Petersson metric.
\endabstract
\thanks
\endthanks
\thanks
\endthanks
\thanks Research  supported by  the National Natural Science Foundation of China (Grant Nos. 11371268, 11631010) and the Natural Science Foundation of Jiangsu Province  (Grant No. BK20141189).
\endthanks
\endtopmatter

\document

\head 1 Introduction and statement of main results\endhead

We begin with some basic definitions and notations.
Let $\Delta=\{z: |z|<1\}$ denote the unit disk in the extended
complex plane $\hat{\Bbb C}$. $\Delta^*=\hat{\Bbb
C}-\overline{\Delta}$ is the exterior of $\Delta$,
$S^1=\partial\Delta=\partial\Delta^*$ is the unit circle, and $\Bbb R$ is the real line. For any function $f=f(z)$ defined on the unit circle $S^1$, we always denote by $\hat f$ the function defined on the real line $\Bbb R$ by $\hat f(\theta)=f(e^{i\theta})$.

Let $\text{Hom}^+(S^1)$ denote the set of all sense-preserving
homeomorphisms of $S^1$ onto itself.
 A homeomorphism
$h\in\text{Hom}^+(S^1)$ is said to be quasisymmetric  if
$$h_{QS}\doteq\sup\{\max(q_h(\theta, t), q^{-1}_h(\theta, t)):\theta\in\Bbb R, t>0\}<+\infty,\tag 1.1 $$
where
$$q_h(\theta, t)\doteq\left|\frac{\hat h(\theta+t)-\hat h(\theta)}{\hat h(\theta)-\hat h(\theta-t)}\right|.\tag 1.2$$
Beurling-Ahlfors [BA] proved that
$h\in\text{Hom}^+(S^1)$ is quasisymmetric if and only if there
exists some quasiconformal homeomorphism  of $\Delta$ onto itself
which has boundary values $h$. Later  Douady-Earle [DE] gave a
quasiconformal extension of $h$ to the unit disk which is
conformally invariant.

The universal Teichm\"uller space $T$ is a universal parameter space
for all Riemann surfaces and can be defined as the right coset space $T=\QS(S^1)/\text{M\"ob}(S^1)$, where $\QS(S^1)$ denotes the group of
all quasisymmetric homeomorphisms of the unit circle, and
$\text{M\"ob}(S^1)$ the subgroup of M\"obius transformations of the
unit disk. The universal Teichm\"uller space $T$ plays
a significant role in Teichm\"uller theory, and it is also a
fundamental object in mathematics and in mathematical physics. On
the other hand, several subclasses of quasisymmetric homeomorphisms
and their Teichm\"uller spaces were introduced and studied for various purposes
in the literature.  We
refer to the books [Ah], [FM], [Ga], [GL], [Hu], [IT], [Le], [Na], [Po2] and the papers [AZ], [Cu], [GS], [FH], [FHS1-2], [HS],  [SW], [TT2], [TWS], [WS] for an introduction to the
subject and more details. In this paper, we are mainly concerned with the so-called Weil-Petersson Teichm\"uller space.

It is well known that  the universal Teichm\"uller space $T$ has a natural complex Banach manifold structure under which the hyperbolic Kobayashi metric is the classical Teichm\"uller metric (see [Ga], [Le], [Na], [Ro]),  and the tangent space to $T$ was identified by Reimann [Re] and later by Gardiner-Sullivan [GS]. Let $\Lambda^*$ denote the Zygmund space in the usual sense (see [Zy]), which consists of all continuous functions $u$ on the unit circle such that
$$
u_{\Lambda^*}\doteq\sup\left\{\frac{|\hat u(\theta+t)-2\hat u(\theta)+\hat u(\theta-t)|}{t}:\theta\in\Bbb R, t>0\right\}<+\infty. \tag 1.3
$$
Then the tangent space to $T$ at the identity map is the set of all functions $u\in \Lambda^*$ which in addition satisfy  the normalized conditions
$$
\Re\bar{\eta}u(\eta)=0,\quad \eta\in S^1
\tag 1.4$$
and
$$
u(1)= u(-1)=u(i)=0.
\tag 1.5$$
More generally, Reimann [Re] proved, given a continuous vector field $u(t, \cdot)\in C^0([0, M], \Lambda^*)$ with the normalized condition (1.4), that the flow maps $h(t, \zeta)$ of the differential equation
$$
\cases
\frac{dh}{dt}=u(t, h)\\
h(0,\zeta)=\zeta
\endcases
\tag 1.6$$
are quasisymmetric homeomorphisms, namely, $h(t, \cdot)\in \QS(S^1)$ for each fixed $t\in [0, M]$.

It is also known that the Kobayashi-Teichm\"uller metric on any Teichm\"uller space  is only induced from a Finsler structure (see [Ob]) and is not a Riemannian metric in general. On the other hand, there does exist a Riemannian metric on a finite dimensional Teichm\"uller space, the Weil-Petersson metric, which has attracted a good bit of attention (see [Hu], [IT], [Mi], [TT2]). In order to extend the definition of the Weil-Petersson metric  to the universal Teichm\"uller space, Nag-Verjovsky [NV] introduced a formal formula for the Weil-Petersson metric, which converges only at those vectors on the unit circle that belong to the Sobolev space $H^{\frac 32}$, however. To overcome this difficulty, Takhtajan-Teo [TT2] endowed the universal Teichm\"uller space with a new complex Hilbert manifold structure, under which the Weil-Petersson metric is a convergent Riemannian metric. But, under this new complex Hilbert manifold structure, the universal Teichm\"uller space $T$ is not connected and has uncountably many connected components. Nowadays, the component containing the identity map is usually called the Weil-Petersson Teichm\"uller space, which is denoted by ${T_0}$ in this paper. Takhtajan-Teo [TT2] proved that, under the Weil-Petersson metric, $T_0$ is precisely the completion of $\Diff_+(S^1)/\text{M\"ob}(S^1)$, the space of all
normalized $C^{\infty}$ diffeomorphisms on the unit circle. Recall that the complex Fr\'echet manifold $\Diff_+(S^1)/\text{M\"ob}(S^1)$ plays an important role in one of the approaches to non-perturbative bosonic closed string field theory based on K\"ahler geometry (see [BR1-2]), and also has an interpretation as a coadjoint orbit of the Bott-Virasoro group (see [Ki], [KY]).

 We say a quasi-symmetric homeomorphism $h$  belongs to the Weil-Petersson class, which is denoted by $\WP(S^1)$, if it represents a point in ${T_0}$. Then ${T_0}=\WP(S^1)/\text{M\"ob}(S^1)$. It is known that  a quasi-symmetric homeomorphism  $h$ belongs to $\WP(S^1)$
 if and only if $h$ has a quasiconformal extension $f$ to
the unit disk whose Beltrami coefficient $\mu$ satisfies the
property that $\iint_{\Delta}|\mu(z)|^2(1-|z|^2)^{-2}dxdy<\infty$. (see [Cu],
[TT2]).
 Due to their importance and wide applications to various areas such as mathematical physics (see [BR1-2], [Ki], [KY], [RSW1-4]),  differential equation and computer vision (see [GMR], [GR], [Ku]), the Weil-Petersson class and its Teichm\"uller space $T_0$ have been much investigated in recent years (see [Fi], [GMR], [GR],  [HS], [Ku], [TT1-2], [Wu]). Recently, motivated by the conformal field theory, Radnell-Schippers-Staubach [RSW1-4] have a programm to extend the Weil-Petersson theory of the universal Teichm\"uller space to the case of Teichm\"uller spaces of bordered Riemann surfaces. Yanagishita (see [Ya1-2] and also [MY]) has even dealt with the Weil-Petersson Teichm\"uller spaces of general Riemann surfaces with a mild geometric condition. However, it is still an open problem how to characterize intrinsically  the elements in $\WP(S^1)$ without using quasiconformal extensions. This problem was proposed by Takhtajan-Teo in 2006 (see page 68 in [TT2]) and was investigated later by  Figalli [Fi], Gay-Balmaz-Marsden-Ratiu ([GMR], [GR]) and some others. In this paper, we will study this problem and  prove the following result, which gives an intrinsic characterization  of a quasisymmetric homeomorphism in the Weil-Petersson class. Recall that, for a function $f$ defined on a set $\Gamma$, $f'$ denotes the derivative of  $f$, namely, for $z\in \Gamma$,
$$f'(z)\doteq\lim_{\Gamma\ni \zeta\to z}\frac{f(\zeta)-f(z)}{\zeta-z} \tag 1.7$$
provided the limit exists, while $f'(z)\doteq 0$ otherwise.

 \proclaim{Theorem 1.1} A sense-preserving  homeomorphism $h$ on the unit circle belongs to the Weil-Petersson class $\WP(S^1)$  if and only if $h$ is absolutely continuous  $($with respect to the arc-length measure$)$ such that $\log h'$ belongs to the Sobolev  class $H^{\frac 12}$.
 \endproclaim

 Theorem 1.1 has several applications which we proceed to explore. It is known that  ${T_0}$ is modeled on the Sobolev  space $H^{\frac 32}$, namely, the tangent space to ${T_0}$ at the identity consists of precisely the $H^{\frac 32}$ vector fields on the unit circle with the normalized conditions (1.2) and (1.3) (see [NV], [TT2]). Recall that when $s>\frac 32$ the group $\Diff^s_+(S^1)$ of all orientation preserving $H^s$ diffeomorphisms of the unit circle and its model space $H^s$ have the same Sobolev $H^s$ regularity.
  An important question is whether the same result holds in the critical case $s=\frac 32$, namely, whether an element in $\WP(S^1)$ also has $H^{\frac 32}$-regularity (see [Fi], [GMR], [GR]). In fact, based on the results by Figalli [Fi], Gay-Balmaz-Marsden-Ratiu ([GMR], [GR])  were able to prove that each homeomorphism in $\WP(S^1)$ belongs to $H^{\frac 32-\epsilon}$ for each $\epsilon>0$. However, we shall prove that the $H^{\frac 32}$-regularity may fail for a quasisymmetric homeomorphism in the Weil-Petersson class, which was conjectured very recently by Gay-Balmaz-Ratiu during their study of the Euler-Weil-Petersson equation (see page 760 in [GR], Conjecture (2)).

 \proclaim{Theorem 1.2} There exists a quasisymmetric homeomorphism  in $\WP(S^1)$ which  belongs neither to the Sobolev class $H^{\frac 32}$ nor to the Lipschitz class $\Lambda^1$.
 \endproclaim

 We will also  deal with the  flows of $H^{\frac 32}$ vector fields on the unit circle.
 It is easy to see that the Weil-Petersson class $\WP(S^1)$ can be generated by the flows of the $H^{\frac 32}$ vector fields on the unit circle (see [GMR], [GR]). However, it is still an open problem whether or not the flows of the $H^{\frac 32}$ vector fields are contained in $\WP(S^1)$, thought it is hoped to be so (see [Fi]). Actually, in the recent paper [GR] by Gay-Balmaz-Ratiu, the authors conjectured that  the flows of the $H^{\frac 32}$ vector fields are contained in $\WP(S^1)$ (see page 760 in [GR], Conjecture and also Conjecture 9.2 below). The following result provides  an affirmative answer to this problem. A more precise statement will be given in Theorem 7.3 below.

 \proclaim {Theorem 1.3} The flows of the $H^{\frac 32}$ vector fields on the unit circle  are always contained in $\WP(S^1)$.
\endproclaim

 Theorem 1.1  is also hoped to be useful to the further study of the geometry and structure of $T_0$. As we shall see later (see Remark 5.1 below), $\WP(S^1)/\Rot(S^1)$ has a very simple model, namely, it can be identified as the real Hilbert space $H^{\frac 12}_{\Bbb R}/{\Bbb R}$ under the bijection  $h\mapsto\log|h'|$.  Here and in what follows, $\Rot(S^1)$ denotes the group of all rotations about the circle $S^1$. Based on this observation, we will  introduce a new metric on $T_0$, which can be defined roughly as follows:
 $$d(h_1, h_2)\doteq\|\log |h'_1|-\log |h'_2|\|_{H^{\frac 12}},\quad h_1, h_2\in {T_0}.\tag 1.8$$
 A precise formula will be given below during the proof of Theorem 1.4 (see (8.5) below). The advantage of this metric is that, as being a global metric, it gives directly the distances between two points in $T_0$. This is in contrast to the case for the Weil-Petersson metric, which is an infinitesimal Riemann metric on the tangent bundle of $T_0$. Anyhow, we shall prove
 \proclaim{Theorem 1.4} The  metric $d$ and the Weil-Petersson metric induce the same topology on ${T_0}$.
 \endproclaim

 We end this Introduction section with the organization of the paper. In Section 2, we give  some basic definitions and results on the universal Teichm\"uller space and the Weil-Petersson Teichm\"uller space. In particular, we establish the complex analytic theory of the pre-logarithmic derivative model of the Weil-Petersson Teichm\"uller space, which plays an important role in the proof of Theorem 1.4. As we shall see later,  our results and their proofs involve much use of the theory of function spaces, in particular, of Sobolev spaces of fractional order. Therefore, in  Sections 3 and 4,   we recall some  basic definitions  on Sobolev spaces,   the BMO space and establish some lemmas that will be frequently used in the proof of  Theorems 1.1 and 1.2; we also deal  with the pull-back operator on $H^{\frac 12}$ by a quasisymmetric homeomorphism and establish several basic results which are needed to prove Theorems 1.1 and 1.3. In Sections 5-8, we give the proofs of Theorems 1.1-1.4. In Section 9, we list several open problems related to this work. In the final Appendix section, we prove Propositions 4.1 and 4.3 stated in Section 4.

 \head 2 Preliminary results on the Weil-Petersson Teichm\"uller space
\endhead

In this section, we give some basic definitions and results on the Weil-Petersson Teichm\"uller space. The results turn out to be essential in the proof of Theorem 1.4. We follow the lines in our recent paper [SW], where the BMO theory of the universal Teichm\"uller space was investigated.

We begin with the standard theory of the universal Teichm\"uller
space (see [Ah], [Le], [Na]). Let  $M(\Delta^*)$ denote  the open
unit ball of the Banach space $L{}^{\infty}(\Delta^*)$ of
essentially bounded measurable functions on $\Delta^*$. For $\mu\in
M(\Delta^*)$, let $f_{\mu}$ be the quasiconformal mapping on the
extended plane $\hat {\Bbb C}$ with complex dilatation equal to
$\mu$ in $\Delta^*$, conformal in $\Delta$, normalized by
$f_{\mu}(0)=f'_{\mu}(0)-1=f''_{\mu}(0)=0$. We say two elements $\mu$
and $\nu$ in $M(\Delta^*)$ are equivalent, denoted by $\mu\sim\nu$,
if $f_{\mu}|_{\Delta}=f_{\nu}|_{\Delta}$. Then
$T=M(\Delta^*)/_{\sim}$ is the Bers model of the universal
Teichm\"uller space. We let $\Phi$ denote the natural  projection
 from $M(\Delta^*)$ onto $T$ so that $\Phi(\mu)$ is the equivalence class $[\mu]$. $[0]$ is called the base point of $T$. The Teichm\"uller distance
between two points $[\mu_1]$ and $[\mu_2]$ in $T$ is defined as
$$\tau([\mu_1], [\mu_2])\doteq\inf\left\{\frac
12\log\frac{1+\|\frac{\nu_1-\nu_2}{1-\overline{\nu}_1\nu_2}\|_{\infty}}
{1-\|\frac{\nu_1-\nu_2}{1-\overline{\nu}_1\nu_2}\|_{\infty}}:
[\nu_1]=[\mu_1], [\nu_2]=[\mu_2] \right\}.\tag 2.1$$

  Let
 $B_2(\Delta)$
denote the Banach space of functions $\phi$ holomorphic in $\Delta$
with norm
$$
\|\phi\|_{B_2}\doteq\sup_{z\in \Delta}(1-|z|^2)^2|\phi(z)|. \tag 2.2
$$
  Consider the map $S: M(\Delta^*)\to B_2(\Delta)$ which sends $\mu$ to
 the Schwarzian derivative of $f_{\mu}|_{\Delta}$. Recall that for any locally univalent function $f$, its Schwarzian derivative $S_f$ is defined by
 $$S_f=N'_f-\frac 12N^2_f,\quad N_f=(\log f')'.\tag 2.3$$
 It is known that
 $S$ is a holomorphic split submersion onto its image, which
descends down to a map $\beta: T\to B_2(\Delta)$  known as the Bers
embedding. Via the Bers embedding, $T$ carries a natural complex Banach manifold
structure so that  $\Phi$  is a holomorphic split submersion.

We proceed to  define the Weil-Petersson Teichm\"uller space (For details, see [TT2] and also [Cu]). We denote by
$\Cal L(\Delta^*)$ the Banach space of all essentially bounded
measurable functions $\mu$ with norm
$$\|\mu\|_{\WP}\doteq\|\mu\|_{\infty}+\left(\frac{1}{\pi}\iint_{\Delta^*}\frac{|\mu(z)|^2}{(|z|^2-1)^2}dxdy\right)^{\frac 12}.\tag 2.4$$
 Set $\Cal
M(\Delta^*)=M(\Delta^*)\cap\Cal L(\Delta^*)$. Then ${T_0}=\Cal
M(\Delta^*)/_{\sim}$ is one of the models of the Weil-Petersson Teichm\"uller space. Actually, ${T_0}$ is the base point component of the universal Teichm\"uller space under the complex Hilbert manifold structure introduced by  Takhtajan-Teo [TT2].

We denote by $\Cal B(\Delta)$ the Banach space of functions $\phi$
holomorphic in $\Delta$ with norm
$$\|\phi\|_{\Cal B}\doteq\left(\frac{1}{\pi}\iint_{\Delta}|\phi(z)|^2(1-|z|^2)^2dxdy\right)^{\frac 12}.\tag 2.5$$
Then,  $\Cal B(\Delta)\subset B_2(\Delta)$, and the
inclusion map is continuous. Under the Bers projection $S: M(\Delta^*)\to
B_2(\Delta)$, $S(\Cal
M(\Delta^*))=S(M(\Delta^*))\cap \Cal B(\Delta)$ (see [Cu], [TT2]). Moreover, we have

\proclaim{Proposition 2.1 ([TT2])} $S:\Cal M(\Delta^*)\to \Cal B(\Delta)$ is a
holomorphic split submersion  from $\Cal M(\Delta^*)$ onto its
image. Consequently, ${T_0}$ has a unique   complex Hilbert manifold structure such that
$\beta: {T_0}\to \Cal B(\Delta)$ is a bi-holomorphic  map from ${T_0}$
onto a domain in $\Cal B(\Delta)$. Under this complex Hilbert manifold structure, the
natural projection $\Phi$ from $\Cal M(\Delta^*)$ onto ${T_0}$ is a
holomorphic split submersion.
\endproclaim

It is well known that a quasiconformal self-mapping of $\Delta^*$
 induces a bi-holomorphic automorphism of the universal
Teichm\"uller space (see [Le], [Na]). Precisely,  let
$w:\Delta^*\to\Delta^*$ be a quasiconformal mapping with complex
dilatation $\mu$. Then $w$ induces an bi-holomorphic isomorphism
$R_w: M(\Delta^*)\to M(\Delta^*)$ as
$$R_w(\nu)=\left(\frac{\nu-\mu}{1-\bar\mu\nu}\frac{\partial w}{\overline{\partial w}}\right)\circ w^{-1}.\tag 2.6$$ $R_w$
descends down a bi-holomorphic isomorphism $w^*: T\to T$ by
$w^*\circ\Phi=\Phi\circ R_w$.

\proclaim{Proposition 2.2} Suppose $w$ is quasi-isometric under the
Poincar\'e metric with Beltrami coefficient  $\mu\in\Cal
M(\Delta^*)$. Then $w^*: {T_0}\to {T_0}$ is bi-holomorphic. \endproclaim
\demo{Proof}Clearly, $R_w$ maps $\Cal M(\Delta^*)$
into itself, and
$R_w:\Cal M(\Delta^*)\to\Cal M(\Delta^*)$ is bi-holomorphic. It
follows from Proposition 2.1  that $w^*: {T_0}\to {T_0}$ is
bi-holomorphic.  \quad$\square$
\enddemo

We continue to consider the  pre-logarithmic derivative model of the Weil-Petersson Teichm\"uller space.
Let
 $B(\Delta)$
denote the  space of functions $\phi$ holomorphic in $\Delta$
with semi-norm
$$
\|\phi\|_{B}\doteq\sup_{z\in \Delta}(1-|z|^2)|\phi'(z)|, \tag 2.7
$$
and $B_0(\Delta)$ the subspace of $B(\Delta)$ which consists of those functions $\phi$ satisfying the condition $\lim_{|z|\to 1}(1-|z|^2)\phi'(z)=0$. Recall that the pre-logarithmic derivative model $\hat T$ of the
universal Teichm\"uller space consists of all functions $\log f'$ (in $B(\Delta)$),
where $f$ belongs to the well known class $S_Q$ of all univalent
analytic functions $f$ in the unit disk $\Delta$ with the normalized
condition $f(0)=f'(0)-1=0$ that can be extended to a quasiconformal
mapping in the whole plane (see [AG], [Zhu]). Under the topology of
Bloch norm (2.7), $\hat T$ is a disconnected open
set. Precisely, $\hat T=\Hat T_b\cup_{\theta\in [0, 2\pi)}\Hat T_{\theta}$, where $\Hat T_b=\{\log f': f\in S_Q\, \text {is
bounded}\}$ and $\Hat T_{\theta}=\{\log f': f\in S_Q\, \text
{satisfies} f(e^{i\theta})=\infty\}$, $\theta\in [0, 2\pi)$, are the
all connected components of $\hat T$ (see [Zhu]). Each $
\Hat T_{\theta}$ is a copy of the Bers model $T$, while $\Hat T_b$ is a
fiber space over $T$. In fact, $\Hat T_b$ is a model of the universal
Teichm\"uller curve (see [Ber], [Te]).

Let $\Cal {AD}(\Delta)$  denote the  space  of all functions $\phi$ holomorphic in $\Delta$ with semi-norm
 $$\|\phi\|_{\Cal {AD}}\doteq\left(\frac{1}{\pi}\iint_{\Delta}|\phi'(z)|^2dxdy\right)^{\frac 12},\tag 2.8$$
 and $\Cal {AD}_0(\Delta)=\{\phi\in\Cal {AD}(\Delta): \phi(0)=0\}$.
 Then $\Cal {AD}(\Delta)\subset B_0(\Delta)$, and the inclusion map is continuous.  We may define $\Cal{AD}(\Delta^*)$  similarly. Define
$$\Lambda(\phi)=\phi''-\frac{1}{2}(\phi')^2,\,\phi\in \Cal {AD}(\Delta).\tag 2.9$$ Then
it holds the following basic result.
\proclaim{Lemma 2.3 ([TT2])} $\Lambda:\Cal {AD}(\Delta)\to\Cal B(\Delta)$ is
holomorphic.
\endproclaim

We come back our situation.  Fix $z_0\in{\Delta^*}$. For $\mu\in M(\Delta^*)$, let
$g^{z_0}_{\mu}$ (abbreviated to  be $g_{\mu}$) be the quasiconformal mapping on the extended plane $\hat
{\Bbb C}$ with complex dilatation equal to $\mu$ in $\Delta^*$,
conformal in $\Delta$, normalized by $g_{\mu}(0)=g'_{\mu}(0)-1=0$,
$g_{\mu}(z_0)=\infty$. Then $\mu\sim\nu$ if and only if $g_{\mu}|_{\Delta}=g_{\nu}|_{\Delta}$. Consider the map $L_{z_0}$ on $ M(\Delta^*)$
by setting $L_{z_0}(\mu)=\log g'_{\mu}$. Then $\cup_{z_0\in\Delta^*}L_{z_0}(M(\Delta^*))=\Hat T_b$, and $\cup_{z_0\in\Delta^*}L_{z_0}(\Cal M(\Delta^*))=\Hat T_b\cap
\Cal {AD}_0(\Delta)$ (see [Cu], [TT2]). We have the following result.

\proclaim{Theorem 2.4} For each $z_0\in{\Delta^*}$,
$L_{z_0}: \Cal M(\Delta^*)\to\Cal {AD}_0(\Delta)$ is holomorphic.
\endproclaim
\demo{Proof} We first  show $L=L_{z_0}: \Cal
M(\Delta^*)\to\Cal {AD}_0(\Delta)$ is continuous.  Recall that $L$ is continuous on $
M(\Delta^*)$ in the topology of Bloch norm (2.7)(see [Le]), namely,
$$\sup_{z\in\Delta}|N_{g_{\nu}}(z)-N_{g_{\mu}}(z)|(1-|z|^2)\le
C(\|\mu\|_{\infty})\|\nu-\mu\|_{\infty},\, \mu, \, \nu\in
M(\Delta^*).\tag 2.10$$
Then,
$$
\align &\|L(\nu)-L(\mu)\|_{\Cal {AD}}^2=\frac{1}{\pi}\iint_{\Delta}|N_{g_{\nu}}(z)-N_{g_{\mu}}(z)|^2dxdy\\
&\le C_1\left(|N_{g_{\nu}}(0)-N_{g_{\mu}}(0)|^2+\iint_{\Delta}(1-|z|^2)^2|N'_{g_{\nu}}(z)-N'_{g_{\mu}}(z)|^2dxdy\right)\\
&\le C_2\left(\|\nu-\mu\|^2_{\infty}+\iint_{\Delta}(1-|z|^2)^2\left|(S_{\nu}(z)-S_{\mu}(z))+\frac{1}{2}(N^2_{g_{\nu}}(z)
-N^2_{g_{\mu}}(z))\right|^2dxdy\right)\\
&\le
C_3(\|\nu-\mu\|^2_{\infty}+\|S(\nu)-S(\mu)\|^2_{\Cal
B}+(\|L(\nu)\|^2_{\Cal {AD}}+\|L(\mu)\|^2_{\Cal {AD}})\|\nu-\mu\|^2_{\infty}).\endalign
$$
By the holomorphy of $S: \Cal M(\Delta^*)\to \Cal B(\Delta)$, we
conclude that $L: \Cal M(\Delta^*)\to\Cal {AD}_0(\Delta)$ is
continuous.

Since $L: \Cal
M(\Delta^*)\to\Cal {AD}_0(\Delta)$ is continuous, we conclude that $L$ is holomorphic by  the infinite dimensional holomorphy (see [Le], [Na]). For completeness, we write down the standard proof. For each $z\in\Delta$, define $l_z(\phi)=\phi(z)$ for $\phi\in
\Cal {AD}_0(\Delta)$. Then, $l_z\in \Cal {AD}^*_0(\Delta)$, that is, $l_z$ is a continuous linear functional on the Banach space $\Cal{AD}_0(\Delta)$. Set $A=\{l_z:
z\in\Delta\}$.  $A$ is a total subset of $\Cal {AD}^*_0(\Delta)$ in the sense that $l_z(\phi)=0$ for all $z\in\Delta$ implies that $\phi=0$.
Now for each $z\in\Delta$, each pair $(\mu, \nu)\in \Cal
M(\Delta^*)\times \Cal L(\Delta^*)$ and small $t$ in the complex
plane,   by the well known holomorphic dependence of quasiconformal
mappings on parameters (see [Ah], [Le], [Na]), we conclude that
$l_z(L(\mu+t\nu))=L(\mu+t\nu)(z)$ is a holomorphic function of $t$.
By a general result about the  infinite
dimensional holomorphy (see [Le], [Na]), it follows that  $L:\Cal M(\Delta^*)\to \Cal {AD}_0(\Delta)$ is holomorphic. \quad$\square$

\enddemo

\proclaim{Theorem 2.5} $\Hat T_b\cap
\Cal {AD}_0(\Delta)$ is a connected open subset of $
\Cal {AD}_0(\Delta)$, and $\Lambda$ is a holomorphic split submersion
from $\Hat T_b\cap\Cal {AD}_0(\Delta)$ onto $\beta({T_0})$.
\endproclaim
\demo{Proof} Clearly,  $\Hat T_b\cap\Cal {AD}_0(\Delta)$ is an
 open subset of $\Cal {AD}_0(\Delta)$. We need  to show that each point of $\Hat T_b\cap
\Cal {AD}_0(\Delta)$ can be connected to 0 by a path in $\Hat T_b\cap
\Cal {AD}_0(\Delta)$.

Let $\log f'\in \Hat T_b\cap
\Cal {AD}_0(\Delta)$. Then $f$ can be extended to a quasiconformal mapping in the whole
plane whose Beltrami coefficient $\mu$ belongs to $\Cal
M(\Delta^*)$, and $z_0=f^{-1}(\infty)\in\Delta^*$. For each $t\in
[0, 1]$,  let $f_t\in S_Q$ be the unique mapping whose
quasiconformal extension to the whole plane has Beltrami coefficient
$t\mu$, and $f_t(z_0)=\infty$. Theorem 2.4 implies that  $\log
f'_t$, $t\in [0, 1]$, is a path in $\Hat T_b\cap
\Cal {AD}_0(\Delta)$
joining $\log f'_0$ to $\log f'$. Now, if $z_0=\infty$, then
$f_0(z)=z$, and we are done. If $z_0\neq\infty$, then
$f_0(z)=z_0z/(z_0-z)$, and $\log f'_0(r\cdot)$, $r\in [0, 1]$, is a
curve in $\Hat T_b\cap
\Cal {AD}_0(\Delta)$  connecting 0 and $\log
f'_0$.

Clearly,  Lemma 2.3 implies that $\Lambda$ is  holomorphic on $\Hat T_b\cap\Cal {AD}_0(\Delta)$. Choose
$z_0\in\Delta^*$. Since $S=\Lambda\circ L_{z_0}$, we conclude that
$\Lambda$ is a holomorphic split submersion
from $\Hat T_b\cap\Cal {AD}_0(\Delta)$ onto $\beta({T_0})$ since $L_{z_0}: \Cal M(\Delta^*)\to
\Hat T_b\cap\Cal {AD}_0(\Delta)$ is holomorphic, and $S: \Cal
M(\Delta^*)\to \Cal \beta({T_0})$ is a holomorphic split submersion. \quad$\square$
\enddemo

\head 3 Some lemmas
\endhead

In this section, we give some lemmas needed to prove Theorems 1.1 and 1.2. First we recall some  basic definitions and results on Sobolev spaces,  the harmonic Dirichlet space and the BMO space that will be frequently used in the rest of the paper (see [Gar], [RS], [Tr]).

For any $s>0$, the Sobolev space $H^s$ consists of all integrable functions $u\in L^1(S^1)$  on the unit circle with semi-norm
$$\|u\|_{H^s}\doteq\left(\sum_{n=-\infty}^{+\infty}|n|^{2s}|a_n(u)|^2\right)^{\frac 12},\tag 3.1$$
where, as usual, $a_n(u)$ is the $n$-th Fourier coefficient of $u$, namely,
$$a_n(u)=\frac{1}{2\pi}\int_0^{2\pi}\hat u(\theta)e^{-in\theta}d\theta. \tag 3.2$$ In this paper, the two cases we are concerned are $s=\frac 32$ and $s=\frac 12$. Recall that $u\in H^{\frac 32}$ if and only if $u$ is absolutely continuous with $u'\in H^{\frac 12}$. It is also known that an integrable function $u$ on the unit circle belongs to $H^{\frac 12}$ if and only if
$$\int_0^{2\pi}\int_0^{2\pi}\frac{|\hat u(s)-\hat u(t)|^2}{|\sin((s-t)/2)|^2}dsdt<+\infty.\tag 3.3$$

We need another description of the space $H^{\frac 12}$. Let $\Cal D(\Delta)$ denote the  space of all harmonic
functions $u$ in the unit disk $\Delta$ with semi-norm
$$\|u\|_{\Cal D}\doteq\left(\frac{1}{\pi}\iint_{\Delta}(|\px u|^2+|\py u|^2)dxdy\right)^{\frac 12}.\tag 3.4$$
 Then, $\Cal
D(\Delta)=\Cal {AD}(\Delta)\oplus\overline{\Cal {AD}(\Delta)}$, or precisely, for each $u\in\Cal
D(\Delta)$, there exists a unique pair of holomorphic  functions $\phi$ and
$\psi$ in $\Cal {AD}(\Delta)$ with $\phi(0)-u(0)={\psi(0)}=0$  such that $u=\phi+\overline\psi$. Here it is a convenient place to introduce two basic operators on the Dirichlet space $\Cal D(\Delta)$. They are $P^+$ and $P^-$, defined respectively by $P^+u=\phi$ and $P^-u=\overline{\psi(\bar z)}$ for $u=\phi+\overline\psi$. It is well known that each
function $u\in\Cal D(\Delta)$ has boundary values almost everywhere on the
unit circle, and the boundary function, still denoted by $u$, belongs to $H^{\frac 12}$, and conversely each function in $H^{\frac 12}$ is obtained in this way  (see [Zy]). In fact, the usual Poisson integral operator $P$ establishes a one-to-one map from $H^{\frac 12}$ onto $\Cal D(\Delta)$ with $\|Pu\|_{\Cal D}=\|u\|_{H^{\frac 12}}$.

Let $I_0$ be a connected (closed) arc on the unit circle $S^1$. An  integrable function $u\in L^1(I_0)$ is said to have bounded mean oscillation if
$$\|u\|_{\BMO(I_0)}\doteq\sup\frac{1}{|I|}\int_I|u(z)-u_I||dz|<+\infty,\tag 3.5$$
where the supremum is taken over all sub-intervals $I$ of $I_0$, while $u_I$ is the average of $u$ on the interval $I$, namely,
$$u_I=\frac{1}{|I|}\int_Iu(z)|dz|.\tag 3.6$$ In particular, $u_{S^1}=a_0(u)$. If $u$ also satisfies the condition
$$\lim_{|I|\to 0}\frac{1}{|I|}\int_I|u(z)-u_I||dz|=0,\tag 3.7$$
we say $u$ has vanishing mean oscillation. These functions are denoted by $\BMO(I_0)$ and $\VMO(I_0)$, respectively. In the following, we are mostly concerned with the case $I_0=S^1$.
Then it is well known that  $H^{\frac 12}\subset \VMO(S^1)$, and the inclusion map is continuous (see [Zh]).

We need some basic results on BMO functions. By the well-known theorem of John-Nirenberg for BMO functions (see [Gar]), there exist two universal positive constants $C_1$ and $C_2$ such that for any $\BMO(I_0)$ function $u$, any subinterval $I$ of $I_0$  and any $\lambda>0$, it holds that
 $$\frac{\left|\{z\in I:|u(z)-u_I|\ge\lambda\}\right|}{|I|}\le C_1\exp\left(\frac{-C_2\lambda}{\|u\|_{\BMO(I_0)}}\right).\tag 3.8$$
 For any $p\ge 1$, by  Chebychev's inequality, we have
 $$
 \align
 \frac{1}{|I|}\int_{I}(e^{|u(z)-u_{I}|}-1)^p|dz|&=\frac{1}{|I|}\int_0^{\infty}\left|\{z\in I:|u(z)-u_{I}|\ge\lambda\}\right|d((e^{\lambda}-1)^p)\\
 &\le pC_1\int_0^{\infty}(e^\lambda-1)^{p-1}e^\lambda\exp\left(\frac{-C_2\lambda}{\|u\|_{\BMO(I_0)}}\right)d\lambda.
 \endalign
 $$
When $p\|u\|_{\BMO(I_0)}<C_2$, we obtain
$$\frac{1}{|I|}\int_{I}(e^{|u(z)-u_{I}|}-1)^p|dz|\le \frac{pC_1\|u\|_{\BMO(I_0)}}{C_2-p\|u\|_{\BMO(I_0)}}.\tag 3.9$$

We will repeatedly use  the following basic result:
\proclaim{Lemma 3.1} Let  $u\in \BMO(I_0)$ and $p\ge 1$.  Then  $e^{u}\in L^p(I_0)$ when $p\|u\|_{\BMO(I_0)}$ is small. In particular, if $u\in \VMO(I_0)$, then  $e^{u}\in L^p(I_0)$ for any real number $p\ge 1$.

\endproclaim
\demo{Proof} When $p\|u\|_{\BMO(I_0)}<C_2$, it follows from (3.9) that
$$\frac{1}{|I_0|}\left\|e^{u-u_{I_0}}-1\right\|^{p}_{p}=\frac{1}{|I_0|}\int_{I_0}|e^{u(z)-u_{I_0}}-1|^p|dz|\le \frac{pC_1\|u\|_{\BMO(I_0)}}{C_2-p\|u\|_{\BMO(I_0)}}.\tag 3.10$$
Consequently,
$$\|e^u\|_p\le e^{\|u\|_1}(\|e^{u-u_{I_0}}-1\|_{p}+|I_0|^{\frac 1p})<+\infty.$$

Now suppose $u\in \VMO(I_0)$, and $p\ge 1$ is any real number. By (3.7), for any sufficiently small subinterval $I$ of $I_0$, $u$ has small BMO norm on $I$ so that $e^u\in L^p(I)$. Decompose $I_0$ as the union of finitely many small subintervals $I_j$ so that $e^u\in L^p(I_j)$, we conclude that $e^u\in L^p(I_0)$ as required.\quad$\square$
 \enddemo

 \proclaim{Lemma 3.2} Let $u\in \VMO(S^1)$ and $u_n\in \BMO(S^1)$ on the unit circle. Suppose that  $\|u_n-u\|_{\BMO(S^1)}\to 0$ and $a_0(u_n-u)\to 0$ when $n\to\infty$, then for any $p\ge 1$, we have
 $\|e^{u_n}-e^u\|_p\to 0$ as $n\to\infty$.
 \endproclaim
 \demo{Proof} By (3.10),
 $$ \left\|e^{(u_n-u)-a_0(u_n-u)}-1\right\|^{2p}_{2p}\le \frac{2pC_1\|u_n-u\|_{\BMO(S^1)}}{C_2-2p\|u_n-u\|_{\BMO(S^1)}}\to 0,\quad n\to\infty.
 $$
  On the other hand, since $u\in \VMO(S^1)$, Lemma 3.1 implies that $e^u\in L^{2p}(S^1)$. Consequently,
 $$
 \align
 \left\|e^{u_n}-e^u\right\|_p
 &\le\left\|e^{u_n-u}-1\right\|_{2p}\|e^u\|_{2p}\\
 &\le\|e^u\|_{2p}\left(e^{a_0(u_n-u)}\|e^{(u_n-u)-a_0(u_n-u)}-1\|_{2p}+\|e^{a_0(u_n-u)}-1\|_{2p}\right),
 \endalign
 $$
 which implies $\|e^{u_n}-e^u\|_p\to 0$ as $n\to\infty$. \quad$\square$
 \enddemo

 Recall that for each
sense-preserving
homeomorphisms $h$ of the unit circle onto itself, there exists some  strictly increasing continuous function $\phi$  on the real line with $\phi(\theta+2\pi)-\phi(\theta)\equiv 2\pi$ such that $h(e^{i\theta})=e^{i\phi(\theta)}$. Then
$$h'(e^{i\theta})=e^{i(\phi(\theta)-\theta)}\phi'(\theta).\tag 3.11$$
Furthermore, $h$ is absolutely continuous on the unit circle if and only if $\phi$ is  absolutely continuous on the real line.

Note that in the statement of Theorem 1.1, the quasi-symmetry of the homeomorphism $h$ is not assumed. The following result gives a sufficient condition for an absolutely continuous sense-preserving  homeomorphism to be quasisymmetric, which will be used in the proof of Theorem 1.1.
\proclaim{Lemma 3.3} Let $h$ be an absolutely continuous sense-preserving  homeomorphism on the unit circle such that  $\log h'\in \VMO(S^1)$. Then $h$ is a quasisymmetric homeomorphism.
\endproclaim
\demo{Proof}  Partyka (see [Pa1], Theorem 3.4.7) asserted that $h$ is actually a symmetric homeomorphism in the sense of Gardiner-Sullivan [GS], namely, for any pair of adjacent subintervals  $I_1$ and $I_2$ in $S^1$ with $|I_1|=|I_2|$, it holds that
$$\frac{|h(I_1)|}{|h(I_2)|}=1+o(1),\quad |I_1|=|I_2|\to 0+.\tag 3.12$$
A detailed proof of this fact was given in [Pa2].  Here we give a fast proof for completeness.

Set $v=\log|h'|$ for simplicity. Then $v\in\VMO(S^1)$.
 For any small subinterval $I$ in $S^1$ such that the BMO-norm of $v$ on $I$ is small, we conclude by (3.9) (with $p=1$)  that
 $$\int_Ie^{|v(z)-v_I|}|dz|\le|I|\left(1+\frac{C_1\|v\|_{\BMO(I)}}{C_2-\|v\|_{\BMO(I)}}\right)=|I|(1+o(1)),\quad |I|\to 0.\tag 3.13$$
 Noting that
$$|h(I)|=\int_I|h'(z)||dz|=\int_Ie^{v(z)}|dz|=e^{v_I}\int_Ie^{v(z)-v_I}|dz|, $$
we obtain from (3.13) that, as $|I|\to 0$,
$$|h(I)|\le e^{v_I}\int_Ie^{|v(z)-v_I|}|dz|\le|I|e^{v_I}(1+o(1)),$$
$$|h(I)|\ge e^{v_I}\int_Ie^{-|v(z)-v_I|}|dz|\ge\frac{|I|^2e^{v_I}}{\int_Ie^{|v(z)-v_I|}|dz|}\ge|I|e^{v_I}(1+o(1)), $$
and so
$$|h(I)|=|I|e^{v_I}(1+o(1)),\quad |I|\to 0.\tag 3.14$$

Now let $I_1$ and $I_2$ be two adjacent subintervals in $[0, 2\pi]$ with $|I_1|=|I_2|=l$ being small such that the BMO-norm of $v$ on $I_1\cup I_2$ is small. It holds that (see [Gar], (1.3) in Chapter VI)
$$|{v_{I_1}}-{v_{I_2}}|=2|{v_{I_1}}-{v_{I_1\cup I_2}}|\le 4\|v\|_{\BMO({I_1\cup I_2})}=o(1),\quad l\to 0+. \tag 3.15$$
Then (3.12)  follows from (3.14-3.15) immediately. $\quad\square $

\enddemo

\proclaim{Lemma 3.4} Let $h$ be an absolutely continuous sense-preserving  homeomorphism on the unit circle. Then $\log h'\in H^{\frac 12}$ if and only if $\log |h'|\in H^{\frac 12}$.
\endproclaim
\demo{Proof} Let $h(e^{i\theta})=e^{i\phi(\theta)}$ as before. Without loss of generality, we assume that $h(1)=1$ so that $\phi(0)=0, \,\phi(2\pi)=2\pi$.  Then $|h'(e^{i\theta})|=\phi'(\theta)$, and
 $$\log h'(e^{i\theta})=\log|h'(e^{i\theta})|+i(\phi(\theta)-\theta).\tag 3.16$$
 It is clear that $\log |h'|\in H^{\frac 12}$ if  $\log h'\in H^{\frac 12}$.

 Conversely, we suppose that $\log |h'|\in H^{\frac 12}$. Set $u=\Im\log h'$ so that $\hat u(\theta)=\phi(\theta)-\theta$. We will show that $u\in H^1$, which implies that $\log h'\in H^{\frac 12}$. In fact, the $n$-th ($n\neq 0$) Fourier coefficient of $u$ is
$$
\align
a_n&=\frac{1}{2\pi}\int_0^{2\pi}\hat u(\theta)e^{-in\theta}d\theta=\frac{1}{2\pi}\int_0^{2\pi}(\phi(\theta)-\theta)e^{-in\theta}d\theta\\
&=\frac{1}{2n\pi i}\int_0^{2\pi}(\phi'(\theta)-1)e^{-in\theta}d\theta=\frac{1}{2n\pi i}\int_0^{2\pi}(|h'(e^{i\theta})|-1)e^{-in\theta}d\theta.
\endalign$$
Thus, by Parseval's equality, we conclude by Lemma 3.1 that
$$\sum_{n\neq 0}n^2|a_n|^2=\frac{1}{4\pi^2}\sum_{n\neq 0}\left|\int_0^{2\pi}(|h'(e^{i\theta})|-1)e^{-in\theta}d\theta\right|^2=\||h'|-1\|^2_2<+\infty.$$
This completes the proof. \quad$\square$
\enddemo

\head 4 Pull-back operator revisited
\endhead

In this section, we deal with the   pull-back operator on the
Sobolev space $H^{\frac 12}$ (and also on the Dirichlet space $\Cal D(\Delta)$)  by a quasisymmetric homeomorphism. The results will be used in the following sections to prove Theorems 1.1 and 1.3 and have independent interests of their own.

Let $h$ be a quasisymmetric homeomorphism.
Then $h$ induces a pull-back operator  by
$$P_hu=u\circ h,\quad u\in H^{\frac 12}.\tag 4.1$$
$P_h$  is a
bounded isomorphism from $H^{\frac 12}$ onto itself with
  $P^{-1}_h=P_{h^{-1}}$. By the well known quasi-invariance of Dirichlet integral under quasiconformal mappings, we have
 $$\|P_h\|\doteq\sup\{\|P_hu\|_{H^{\frac 12}}: \|u\|_{H^{\frac 12}}=1\}\le e^{\tau(0, h)}.\tag 4.2$$
This operator has played an important role in the study of Teichm\"uller theory (see [HS], [NS], [Pa1], [SW], [TT2]). As stated in the introduction, the universal Teichm\"uller space has a quasisymmetric homeomorphism model, namely, $T=\QS(S^1)/\text{M\"ob}(S^1)$.   Nag-Sullivan [NS]
proved that   the universal Teichm\"uller space $T$ can be embedded in the
universal Siegel period matrix space by means of the operator $P_h$ (see also [TT2]). Notice that $P_h$ (or more precisely, $P\circ P_h$, the composition of $P_h$ with the Poisson integral operator $P$) is also a
bounded isomorphism from $\Cal D(\Delta)$ onto itself, and $P^{-1}_h=P_{h^{-1}}$.

We will need the following result.
\proclaim{Proposition 4.1} Let $h$ and $h_0$ be quasisymmetric homeomorphisms which keep the points $1$, $-1$ and $i$ fixed. Then for each fixed $u\in H^{\frac 12}$, $\|P_hu-P_{h_0}u\|_{H^{\frac 12}}\to 0$ when $\tau(h, h_0)\to 0$.
\endproclaim
As far as the author know,  Proposition 4.1 is not available in the literature. We will prove it in the final Appendix section. A natural question to ask is
\proclaim{Question 4.2} Under the assumption of Proposition 4.1, is it true that $\|P_h-P_{h_0}\|\to 0$ when $\tau(h, h_0)\to 0$?
\endproclaim
We proceed  to investigate   the pull back operator $P_h$ induced by a quasisymmetric homeomorphism.
When restricted to $\Cal {AD}(\Delta)$, $P_h$ (more precisely, $P\circ P_h$) is a bounded operator from $\Cal
{AD}(\Delta)$ into $\Cal D(\Delta)$. So we may define two further operators
$P^+_h=P^+\circ P_h$ and $P^-_h=P^-\circ P_h$. Both $P^+_h$ and
$P^-_h$ are bounded operators from $\Cal {AD}(\Delta)$ into itself. For completeness, we recall that  $P^-_h$ is a compact operator if and only if
$h$ is symmetric, while $P^-_h$ is a Hilbert-Schmidt operator if and
only if $h$ belongs to the Weil-Petersson class $\WP(S^1)$ (see [HS]).  We will not use this result in this paper.

The following result will play an important role in the proof of Theorem 1.1.

\proclaim {Proposition 4.3} $P^+_h$ is a bounded isomorphism from $\Cal{AD}(\Delta)$ onto itself.  Moreover, it holds that
$$\|P^+_h\phi\|^2_{\Cal{AD}}=\|\phi\|^2_{\Cal{AD}}+\|P^-_h\phi\|^2_{\Cal{AD}},\,\phi\in\Cal{AD}(\Delta).\tag 4.3$$
\endproclaim

Proposition 4.3 may be a known result, but, to the best of the author's knowledge, a proof does not
appear in the literature. We will give the proof in the final Appendix section.

We now establish a technical result used to prove Theorem 1.1. We consider the harmonic conjugation operator $H$ in the usual sense. Precisely, for a real valued integrable  function $u$ on the unit circle, there exists a unique harmonic function $v$ on the unit disk with $v(0)=0$ such that $Pu+iv$ is analytic. Then $Hu=v|_{S^1}$. When $u$ is complex valued, set $Hu=H\Re u+iH\Im u$. Then, $\overline{Hu}=H\overline u$, and
$H\phi=-i(\phi-\phi(0))$ when $\phi$ is holomorphic.
We have the following basic result:
\proclaim{Lemma 4.4} For each $\phi\in\Cal{AD}(\Delta)$, it holds that
$$(HP_h+P_hH)\phi=-i(2P^+_h\phi-P^+_h\phi(0)-\phi(0)).$$
\endproclaim
\demo{Proof} The proof goes  as follows:
$$
\align
(HP_h+P_hH)\phi(z)&=H(P^+_h\phi(z)+{P^-_h\phi(\bar z)})-iP_h(\phi(z)-\phi(0))\\
&=-i(P^+_h\phi(z)-P^+_h\phi(0))+i{P^-_h\phi(\bar z)}-i(P^+_h\phi(z)+{P^-_h\phi(\bar z)}-\phi(0))\\
&=-i(2P^+_h\phi(z)-P^+_h\phi(0)-\phi(0)). \quad\square
\endalign
$$

\enddemo

\proclaim{Corollary 4.5} Let $v\in H^{\frac 12}$ be real valued.  Then there exists some $u\in H^{\frac 12}$ such that $\|(HP_h+P_hH)u-v\|_{H^{\frac 12}}=0$. Furthermore, $2\|u\|_{H^{\frac 12}}\le \|v\|_{H^{\frac 12}}$.
\endproclaim

\demo{Proof} Set $\psi=i(v+iHv)/2$. Then $P\psi\in \Cal{AD}(\Delta)$. By Proposition 4.3, there exists $\phi\in\Cal{AD}(\Delta)$ such that $P^+_h\phi=P\psi$. Letting $u=\Re\phi$, we obtain by Lemma 4.4 that
$$(HP_h+P_hH)u=\Re(HP_h+P_hH)\phi=\Im(2P^+_h\phi-(P^+_h\phi(0)-\phi(0)))=v-\Im(P^+_h\phi(0)-\phi(0)).$$
Consequently, $\|(HP_h+P_hH)u-v\|_{H^{\frac 12}}=0$, and  by (4.3),
$$4\|u\|^2_{H^{\frac 12}}=2\|\phi\|^2_{\Cal{AD}} \le 2\|P\psi\|^2_{\Cal{AD}}=\|v\|^2_{H^{\frac 12}}. \quad\square $$

\enddemo

\head 5 Proof of Theorem 1.1
\endhead

In this section, we will give the proof of  Theorem 1.1. We first recall the normalized decomposition of a quasisymmetric homeomorphism. For any quasisymmetric homeomorphism $h$, there exists a unique pair of conformal mappings $f\in S_Q$ and $g$ on $\Delta$ and $\Delta^*$, respectively, such that $f(0)=f'(0)-1=0$, $g(\infty)=\infty$,  $h=f^{-1}\circ g$ on $S^1$. We call this a normalized decomposition of $h$. Conversely, for each $f\in S_Q$ which maps the unit disk onto a bounded Jordan domain, there exists a quasisymmetric $h$ with the normalized decomposition $h=f^{-1}\circ g$.  It is clear that $h$ is uniquely determined if $h(1)=1$, and in this case we say  $h$ is the normalized conformal sewing mapping of $f$.

\vskip 0.3 cm
\noindent {\bf {Proof of \lq\lq only if \rq\rq part:}} Suppose $h\in \WP(S^1)$. Consider the above normalized decomposition $h=f^{-1}\circ g$. Then, $\log f'\in\Cal{AD}(\Delta)$,  $\log g'\in\Cal {AD}(\Delta^*)$. For details, see [TT2] and also [Cu]. Then, $h$ is absolutely continuous on $S^1$, and from $f\circ h=g$ we obtain $(f'\circ h)h'=g'$. Thus,
$$\log h'=\log g'-\log f'\circ h=\log g'-P_h\log f'.\tag 5.1$$
Consequently, $\log h'\in H^{\frac 12}$. \quad$\square$

\vskip 0.3 cm

\noindent {\bf {Proof of \lq\lq  if\rq\rq part:}}  The proof of this direction is more difficult. Suppose $h$ is an absolutely continuous homeomorphism on the unit circle such that $\log h'\in H^{\frac 12}$. Lemma 3.3 implies that $h$ is a quasisymmetric homeomorphism so that Corollary 4.5 may be used.  Without loss of generality, we assume $h(1)=1$. Then $h(e^{i\theta})=e^{i\phi(\theta)}$, where $\phi$ is a strictly  increasing and   absolutely continuous function
on the real line $\Bbb R$  such that $\phi(0)=0, \, \phi(\theta+2\pi)-\phi(\theta)\equiv 2\pi.$

We first assume $\|\log h'\|_{H^{\frac 12}}$ is small. By Corollary 4.5, there exists
 some $u\in H^{\frac 12}$ and a real constant $c_1$ such that
 $$(HP_h+P_hH)u=-H\log |h'|-\Im\log h'+c_1,\tag 5.2$$
and
 $2\|u\|_{H^{\frac 12}}\le \|H\log |h'|+\Im\log h'\|_{H^{\frac 12}}$ is small.
 Then there exists a locally univalent  analytic function $f$ on the unit disk with $f(0)=f'(0)-1=0$ such that for some constant $c_2$,
 $$\log f'(z)=P(u+iHu)(z)+c_2.\tag 5.3$$
 Since $\|\log f'\|_{\Cal {AD}}=\|u+iHu\|_{H^{\frac 12}}$ is small, by the continuity of the inclusion of $\Cal{AD}(\Delta)$ into $B(\Delta)$, $\|\log f'\|_B$ is also small. It is well known that $f$ is univalent in $\Delta$ and can be extended to a quasiconformal mapping in the whole plane (see [Be] and also [AG]). Consequently, $\log f'\in \Hat T_b\cap\Cal{AD}_0(\Delta)$.

 Now we set $v=P_hu+\log |h'|$. Then $\|v\|_{H^{\frac 12}}$ is  small. In fact, when $\|\log h'\|_{H^{\frac 12}}$ is small, $\|\log h'\|_{\BMO(S^1)}$ is also small by the continuity of the inclusion $H^{\frac 12}$ into $\VMO(S^1)\subset\BMO(S^1)$. Then $h$ can be extended to a quasiconformal mapping in the unit disk whose Beltrami coefficient $\mu$ has small norm $\|\mu\|_{\infty}$ (see [AZ], [Be]), which in turn implies by (4.2) that $\|P_hu\|_{H^{\frac 12}}$ is small and so $\|v\|_{H^{\frac 12}}$ is also small. By the same reasoning as above, there exists a quasiconformal mapping $g$ on the whole plane with $g(\infty)=\infty$ such that $g$  is conformal in $\Delta^*$ with $\log g'\in\Cal {AD}(\Delta^*)$ and
 $$\log g'=v-iHv+(c_2+ic_1)=P_hu+\log |h'|-iHP_hu-iH\log |h'|+c_2+ic_1.\tag 5.4$$

 Now it follows from (5.2-5.4) that
 $$\align
 P_h\log f'-\log g'&=(P_hu+iP_hHu+c_2)-(P_hu+\log |h'|-iHP_hu-iH\log |h'|+c_2+ic_1)\\
 &=i(P_hHu+HP_hu)-\log |h'|+iH\log |h'|-ic_1\\
 &=-i(H\log |h'|+\Im\log h')-\log |h'|+iH\log |h'|\\
 &=-\log h'.
 \endalign
 $$
 Consequently, adding some constant to $g$ if necessary, it holds that $g=f\circ h$. Since  $\log f'\in \Hat T_b\cap\Cal{AD}_0(\Delta)$, we conclude that $h$ belongs to the Weil-Petersson class under the assumption that $\|\log h'\|_{H^{\frac 12}}$ is small. It should be pointed out that the above reasoning was inspired by David [Da] in an other  setting of BMO theory of the universal Teichm\"uller space.

 When $\|\log h'\|_{H^{\frac 12}}$ is not necessarily small, we use an approximation process. Since  $\log h'\in H^{\frac 12}$, there exists a sequence $(u_n)$ of real valued (real) analytic  functions such that $\|u_n-\log|h'|\|_{H^{\frac 12}}\to 0$ as $n\to\infty$. Replacing $u_n$ by $u_n-a_0(u_n)+a_0(\log|h'|)$ if necessary, we may assume that $a_0(u_n)=a_0(\log|h'|)$. Define $h_n(e^{i\theta})=e^{i\phi_n(\theta)}$ by
 $$\phi_n(\theta)=\frac{2\pi}{\int_0^{2\pi}e^{\hat{u_n}(t)}dt}\int_0^{\theta}e^{\hat{u_n}(t)}dt,\, \theta\in \Bbb R.\tag 5.5$$
   Then, $h_n\in \WP(S^1)$ since $\phi_n$ is a real analytic diffeomorphism.

 We first show that $\|\log h_n'-\log h'\|_{H^{\frac 12}}\to 0$ as $n\to\infty$. By our construction, $\|\log|h_n'|-\log|h'|\|_{H^{\frac 12}}\to 0$ as $n\to\infty$. We need to show that $\|\Im\log h_n'-\Im\log h'\|_{H^{\frac 12}}\to 0$ as $n\to\infty$. For simplicity, we set $\lambda_n=\Im\log h_n'-\Im\log h'$  so that $\hat{\lambda_n}=\phi_n-\phi$. Recall that $H^{\frac 12}\subset \VMO(S^1)$, and the inclusion map is continuous. Noting that
 $$|a_0(e^{u_n})-1|=\frac{1}{2\pi}\left|\int_0^{2\pi}(e^{\hat{u_n}(t)}-e^{\log\phi'(t)})dt\right
 |\le \|e^{u_n}-e^{\log|h'|}\|_1, $$
we conclude by Lemma 3.2 that $a_0(e^{u_n})\to 1$ as $n\to\infty$. Now (5.5) implies that
$\log|h'_n|=u_n-\log a_0(e^{u_n})$, which implies $a_0(\log|h'_n|)=a_0(u_n)-\log a_0(e^{u_n})\to a_0(\log|h'|)$ as $n\to\infty$.
By Lemma 3.2 again, we conclude that, for any $p\ge 1$, $\||h_n'|-|h'|\|_p\to 0$ as $n\to\infty$. Now the $m-$th $(m\neq 0)$ Fourier coefficient of $\lambda_n$ is
$$
\align
a_m&=\frac{1}{2\pi}\int_0^{2\pi}\hat{\lambda_n}(\theta)e^{-im\theta}d\theta=\frac{1}{2\pi}\int_0^{2\pi}(\phi_n(\theta)-\phi(\theta))e^{-im\theta}d\theta\\
&=\frac{1}{2m\pi i}\int_0^{2\pi}(\phi'_n(\theta)-\phi'(\theta))e^{-im\theta}d\theta=\frac{1}{2m\pi i}\int_0^{2\pi}(|h'_n|(e^{i\theta})-|h'|(e^{i\theta}))e^{-im\theta}d\theta,
\endalign$$
we conclude by Parseval's equality that
$$\|\lambda_n\|_{H^{1}}=\sum_{m\neq 0}m^2|a_m|^2=\frac{1}{4\pi^2}\sum_{m\neq 0}\left|\int_0^{2\pi}(|h'_n|(e^{i\theta})-|h'|(e^{i\theta}))e^{-im\theta}d\theta\right|^2=\||h'_n|-
|h'|\|^2_2, $$
which implies $\|\lambda_n\|_{H^{\frac 12}}\le \|\lambda_n\|_{H^{1}}\to 0$ as $n\to\infty$. Thus,  $\|\log h_n'-\log h'\|_{H^{\frac 12}}\to 0$ as $n\to\infty$.

Now we consider $\tilde h_n=h_n\circ h^{-1}$. Then  $\tilde h_n$ is absolutely continuous. Noting that
$$\log\tilde h'_n=(\log h'_n-\log h')\circ h^{-1}=P^{-1}_h(\log h'_n-\log h'),$$
we find that $\|\log\tilde h_n'\|_{H^{\frac 12}}\to 0$ as $n\to\infty$. By what we have proved in the small norm case, $\tilde h_n\in \WP(S^1)$. Since $\WP(S^1)$ is a group (see [Cu], [TT2]), we conclude that $h\in \WP(S^1)$. Now the proof of Theorem 1.1 is completed. \quad$\square$

\vskip 0.3 cm
\noindent {\bf Remark 5.1:}\quad By means of Theorem 1.1, we can give a new model of the Weil-Petersson Teichm\"uller space. More precisely, let $H^{\frac 12}_{\Bbb R}$ denote the subspace of all real-valued functions in $H^{\frac 12}$. By Theorem 1.1, $\log |h'|\in H^{\frac 12}_{\Bbb R}$ for $h\in \WP(S^1)$. Conversely, suppose $u\in H^{\frac 12}_{\Bbb R}$. Adding to a constant if necessary, we may assume that $\int_0^{2\pi}e^{\hat u(t)}dt=2\pi$. Set $h(e^{i\theta})=e^{i\phi(\theta)}$ by
$$\phi(\theta)=\int_0^{\theta}e^{\hat u(t)}dt,\, \theta\in \Bbb R.\tag 5.6$$
Then $h$ is an absolutely continuous sense-preserving homeomorphism of the unit circle with $\log|h'|=u$. By Lemma 3.4 and Theorem 1.1, we get $h\in \WP(S^1)$. Consequently, the correspondence $h\mapsto \log|h'|$ establishes a one-to-one map from $\WP(S^1)/\Rot(S^1)$ onto $H^{\frac 12}_{\Bbb R}/{\Bbb R}$. By means of the $H^{\frac 12}$ metric, a  metric can be assigned to $\WP(S^1)/\Rot(S^1)$. This will be done in Section 8.

\head 6 A counterexample: Proof of Theorem 1.2
\endhead

 Combining with  Theorem 1.1, the following result gives the proof of  Theorem 1.2.
\proclaim{Theorem 6.1} Fix $\alpha>1$. Define  $h(e^{i\theta})=e^{i\varphi(\theta)}$
$$\varphi(\theta)=c_{\alpha}\int_0^{\theta}\left(\left(\log \alpha-\log\sin\frac{t}{2}\right)^2+\frac{(\pi-t)^2}{4}\right)dt,\,\theta\in [0, 2\pi], \tag 6.1$$
 where $c_{\alpha}>0$ is a constant so that $\varphi(2\pi)=2\pi$. Then $h$ is a sense-preserving  homeomorphism  which is absolutely continuous such that $\log h'\in H^{\frac 12}$, but $h$ is neither $H^{\frac 32}$ nor Lipschitz.
\endproclaim
\demo{Proof} We first point out that $\varphi$ can be extended to the whole real line  $\Bbb R$ by means of  $\varphi(\theta+2\pi)-\varphi(\theta)\equiv 2\pi$. Consider
$$g(z)=\log\log\frac{2\alpha}{1-z}.\tag 6.2$$
$g$ is holomorphic in $\Delta$, and except for $e^{i\theta}=1$, $\lim_{z\to e^{i\theta}}g(z)$ exists and equals
$$g(e^{i\theta})=\log\left(\log \alpha-\log\sin\frac{\theta}{2}+i\frac{\pi-\theta}{2}\right).$$
We first show that $g\in\Cal {AD}(\Delta)$. Noting that
$$g'(z)=\frac{1}{(1-z)\log\frac{2\alpha}{1-z}},$$
it is sufficient to show that $$\iint_{\{|z-1|<1\}}|g'(z)|^2dxdy<+\infty.$$
This can be be done as follows:
$$
\align
\iint_{\{|z-1|<1\}}|g'(z)|^2dxdy&=\int_{\{|w|<1\}}\frac{1}{|w\log\frac{2\alpha}{ w}|^2}dudv\\
&=\int_0^1\rho d\rho\int_0^{2\pi}\frac{d\theta}{\rho^2(\log^2\frac{2\alpha}{\rho}+\theta^2)}\\
&=\int_0^1\frac{1}{\rho\log\frac{2\alpha}{\rho}}\arctan\frac{2\pi}{\log\frac{2\alpha}{\rho}}d\rho\\
&=\int_{\log 2\alpha}^{+\infty}\frac{\arctan\frac{2\pi}{x}}{x}dx\\
&<2\pi\int_{\log 2\alpha}^{+\infty}\frac{1}{x^2}dx=\frac{2\pi}{\log 2\alpha}.
\endalign
$$
Thus,  $g\in H^{\frac 12}$, which implies that  ${\Re}\,g\in H^{\frac 12}$. By Lemma 3.1, we obtain that $\exp(2\Re g)\in L^1(S^1)$. Noting that
$${\Re}\,g(e^{i\theta})=\log\left|\log \alpha-\log\sin\frac{\theta}{2}+i\frac{\pi-\theta}{2}\right|=\frac 12\log\left(\left(\log \alpha-\log\sin\frac{\theta}{2}\right)^2+\frac{(\pi-\theta)^2}{4}\right)$$
when $\theta\in (0, 2\pi)$, we conclude that our function $\varphi$ defined in (6.1) is well-defined, strictly increasing  and absolutely continuous with
$$\varphi'(\theta)=c_{\alpha}\left(\left(\log \alpha-\log\sin\frac{\theta}{2}\right)^2+\frac{(\pi-\theta)^2}{4}\right),\,\theta\in (0, 2\pi).\tag 6.3$$
Thus, $h$ is an absolutely continuous sense-preserving homeomorphism of the unit circle onto itself. Since  $\|\varphi'\|_{\infty}=\infty$,  $h$ is not Lipschitz.
On the other hand, since $\log|h'|=\log c_{\alpha}+2\Re g\in H^{\frac 12}$, we conclude by Lemma 3.4 that $\log h'\in H^{\frac 12}$.

It remains to  show that $h$ is not in $H^{\frac 32}$, or equivalently, $h'$ is not in $H^{\frac 12}$. By means of  (3.3), it is sufficient to show that $|h'|$ is not in $H^{\frac 12}$. To do so,  we consider the following analytic function in the unit disk
$$f(z)=\log(1-z).\tag 6.4$$ Then, except for $e^{i\theta}=1$, $\lim_{z\to e^{i\theta}}f(z)$ exists and is equal to
$$f(e^{i\theta})=\log(1-e^{i\theta})=\log 2+\log\sin\frac{\theta}{2}-i\frac{\pi-\theta}{2}.\tag 6.5$$
It is easy to see that $f$ does not belong to $\Cal{AD}(\Delta)$, which implies that ${{\Re}}\, f$ is not in $H^{\frac 12}$. By (3.3) we have
$$\int_0^{\pi}\int_0^{\pi}\frac{|\log\sin s-\log\sin t|^2}{|\sin(s-t)|^2}dsdt=+\infty.\tag 6.6$$

Fix $0<\epsilon<\pi/4$, and set $I_{\epsilon}=[\pi/2-\epsilon, \pi/2+\epsilon]$, $(I_{\epsilon}\times I_{\epsilon})^c=[0, \pi]\times[0,
\pi]-I_{\epsilon}\times I_{\epsilon}$.  Noting that $\log(1+x)<x$ when $x>0$, we find that
$$|\log x-\log y|\le\frac{|x-y|}{\min(x, y)},\, x>0, y>0.$$
On the other hand, $\sin x\ge (2/\pi) x$ when $0<x<\pi/2$, we conclude that
$$\frac{|\log\sin s-\log\sin t|^2}{|\sin(s-t)|^2}\le\frac{\pi^2}{4}\frac{|\sin s-\sin t|^2}{|s-t|^2\min(\sin^2s, \sin^2t)}\le\frac{\pi^2}{4\cos^2\epsilon}\le\frac{\pi^2}{2}, \, s,\,t\in I_{\epsilon}.$$
Thus,
$$\int_{I_{\epsilon}}\int_{I_{\epsilon}}\frac{|\log\sin s-\log\sin t|^2}{|\sin(s-t)|^2}dsdt<+\infty.$$
It follows from (6.6) that
$$\iint_{(I_{\epsilon}\times I_{\epsilon})^c}\frac{|\log\sin s-\log\sin t|^2}{|\sin(s-t)|^2}dsdt=+\infty.$$
Noting that $\log\sin s<\log\cos\epsilon<0$ when $s\in I^c_{\epsilon}$, we conclude  from the above equality  that
$$\align
&\iint_{(I_{\epsilon}\times I_{\epsilon})^c}\frac{|(\log \alpha-\log\sin s)^2-(\log \alpha-\log\sin t)^2|^2}{|\sin(s-t)|^2}dsdt\\
&\ge \log^2(\alpha^2\cos\epsilon)\iint_{(I_{\epsilon}\times I_{\epsilon})^c}\frac{|\log\sin s-\log\sin t|^2}{|\sin(s-t)|^2}dsdt=+\infty,
\endalign$$
which implies that
$$\int_0^{\pi}\int_0^{\pi}\frac{|(\log \alpha-\log\sin s)^2-(\log \alpha-\log\sin t)^2|^2}{|\sin(s-t)|^2}dsdt=+\infty.\tag 6.7$$

On the other hand, consider the function $u$ on the unit circle defined by $u(e^{i\theta})=(\pi-\theta)^2$, $\theta\in [0, 2\pi]$. Then, $u\in H^{\frac 12}$. Actually, a direct computation will show that the $n$-th ($n\neq 0$) Fourier coefficient of $u$ is
$$a_n=\frac{1}{2\pi}\int_0^{2\pi}\hat u(\theta)e^{-in\theta}d\theta=\frac{1}{2\pi}\int_0^{2\pi}(\pi-\theta)^2e^{-in\theta}d\theta=\frac {2}{n^2}.$$ Combining this with (3.3) and (6.7), we conclude that  $|h'|$ is not in $H^{\frac 12}$.
This completes the proof of Theorem 6.1. \quad$\square$

\enddemo

\head 7 Proof of Theorem 1.3
\endhead

We first prove two general results.

\proclaim{Lemma 7.1} Given a continuous vector field $u(t, \cdot)\in C^0([0, M], \Lambda^*)$ with the normalized conditions (1.4) and (1.5),  the flow maps $h(t, \zeta)$ of the differential equation
$$
\cases
\frac{dh}{dt}=u(t, h)\\
h(0,\zeta)=\zeta
\endcases
\tag 7.1$$
are quasisymmetric homeomorphisms, and $h(t, \cdot): [0, M]\rightarrow T$ is continuous.
\endproclaim

\demo{Proof} As stated in Section 1,  Reimann [Re] proved that, for each fixed $t\in [0, M]$, $h(t, \cdot)$ is a quasisymmetric homeomorphism. In fact, Agard-Kelingos (see Theorems 1 and 2 in [AK]) already proved that $h(t, \cdot): [0, M]\rightarrow T$ is continuous under the assumption that $u(t, \cdot)$ can be extended to a so-called quasiconformal deformation $U(t, \cdot)$ to the unit disk with $\overline{\partial}U(t, \cdot)\in C^0([0, M], L^{\infty}(\Delta))$, which was proved to be true by Gardiner-Sullivan (see Section 8 in [GS])  and Reich-Chen (see Theorem 2.2 in [RC]) independently. A detailed proof of Lemma 7.1 can be found  in our paper [HWS].\quad$\square$
\enddemo

\proclaim{Lemma 7.2} Let $h_t$, $t\in [0, M]$ be quasisymmetric homeomorphisms which keep the points $1$, $-1$ and $i$ fixed. Suppose $u_t:[0, M]\to H^{\frac12}$ and $h_t:[0, M]\to T$ are continuous. Then $P_{h_t}u_t: [0, M]\to H^{\frac 12}$ is continuous.
\endproclaim
\demo{Proof} Fix $t_0\in [0, M]$. By (4.2)  we have
$$
\align
\|P_{h_t}u_t-P_{h_{t_0}}u_{t_0}\|_{H^{\frac 12}}&\le \|P_{h_t}u_t-P_{h_{t}}u_{t_0}\|_{H^{\frac 12}}+\|P_{h_t}u_{t_0}-P_{h_{t_0}}u_{t_0}\|_{H^{\frac 12}}\\
&\le e^{\tau(0, h_t)}\|u_t-u_{t_0}\|_{H^{\frac 12}}+\|P_{h_t}u_{t_0}-P_{h_{t_0}}u_{t_0}\|_{H^{\frac 12}}.
\endalign
$$
We conclude that $P_{h_t}u_t: [0, M]\to H^{\frac 12}$ is continuous by Proposition 4.1 and the continuity of $u_t$ and $h_t$. \quad$\square$

\enddemo
Now we begin to prove Theorem 1.3. It is contained in
 \proclaim {Theorem 7.3} Given a continuous vector field $u(t, \cdot)\in C^0([0, M], H^{\frac 32})$ with the normalized condition $(1.4)$,  the flow maps $h(t, \cdot)$ of the differential equation $(7.1)$
belong to the Weil-Petersson class, namely, $h(t, \cdot)\in \WP(S^1)$ for each fixed $t\in [0, M]$; Furthermore, the mapping $t\mapsto \log h'(t, \cdot)$ from $[0, M]$ into $H^{\frac 12}$ is continuously differentiable such that
$$\frac{d}{dt}\log h'(t, \cdot)=u'(t, h(t, \cdot)).\tag 7.2$$
\endproclaim
\demo {Proof} Without loss of the generality, we assume that the vector field $u(t, \cdot)$ also satisfies the normalized condition (1.5) so that the the flow maps $h(t, \cdot)$ keep the points $1$, $-1$ and $i$ fixed. We first point out that by Figalli's result (see [Fi]), for each fixed $t\in [0, M]$, $h(t, \cdot)$ is absolutely continuous. As done by Figalli [Fi], differentiating both sides of the equation
$$\frac{d}{dt}h(t, \zeta)=u(t, h(t, \zeta))\tag 7.3$$
with respect to $\zeta$ yields
$$\frac{d}{dt}h'(t, \zeta)=u'(t, h(t, \zeta))h'(t, \zeta),$$
that is,
$$\frac{d}{dt}\log h'(t, \zeta)=u'(t, h(t, \zeta)).$$
Noting that $h(0, \zeta)=\zeta$, we obtain
$$\log h'(t, \zeta)=\int_0^tu'(s, h(s, \zeta))ds.\tag 7.4$$
Recalling that the inclusion of $H^{\frac 32}$ into $\Lambda^*$ is continuous, we conclude  that $u'(t, h(t, \cdot)): [0, M]\rightarrow H^{\frac 12}$ is continuous by Lemmas 7.1 and 7.2.  Now Theorem 7.3  follows from the following Lemma 7.4 immediately.
\quad$\square$

\enddemo

\proclaim {Lemma 7.4} Suppose $u(t, \cdot): [0, M]\rightarrow H^{\frac 12}$ is continuous, and
$$U(t, \zeta)=\int_0^tu(s, \zeta)ds,\quad \zeta\in S^1.\tag 7.5$$
Then for each fixed $t\in [0, M]$, $U(t, \cdot)\in H^{\frac 12}$, and $U(t, \cdot): [0, M]\rightarrow H^{\frac 12}$ is continuous differentiable with
$$\frac {d}{dt}U(t, \cdot)=u(t, \cdot).\tag 7.6$$
\endproclaim
\demo{Proof} For simplicity, we set $U(t, \cdot)=U_t$, $u(t, \cdot)=u_t$. By definition we have
$$
\align
a_n(U_t)&=\frac{1}{2\pi}\int_0^{2\pi}\left(\int_0^tu_s(e^{i\theta})ds \right)e^{-in\theta}d\theta\\
&=\int_0^t\left(\frac{1}{2\pi}\int_0^{2\pi}u_s(e^{i\theta})e^{-in\theta}d\theta \right)ds\\
&=\int_0^ta_n(u_s)ds.
\endalign$$
Then,
$$|a_n(U_t)|^2\le t\int_0^t|a_n(u_s)|^2ds,$$
$$\align
\|U_t\|^2_{H^{\frac 12}}&=\sum_{n=-\infty}^{+\infty}|n||a_n(U_t)|^2\\
&\le t\sum_{n=-\infty}^{+\infty}|n|\int_0^t|a_n(u_s)|^2ds\\
&=t\int_0^t\sum_{n=-\infty}^{+\infty}|n||a_n(u_s)|^2ds\\
&=t\int_0^t\|u_s\|^2_{H^{\frac 12}}ds\\
&\le t^2\max_{s\in [0, t]}\|u_s\|^2_{H^{\frac 12}}.
\endalign
$$
Consequently, for each fixed $t\in [0, M]$, $U(t, \cdot)\in H^{\frac 12}$.

It remains to show (7.6). Fix $t_0\in [0, M]$. Noting that
$$U_{t_0+t}(\zeta)-U_{t_0}(\zeta)-tu_{t_0}(\zeta)=\int_{t_0}^{t_0+t}(u_s(\zeta)-u_{t_0}(\zeta))ds,$$
we conclude by the reasoning as above that
$$\|U_{t_0+t}-U_{t_0}-tu_{t_0}\|_{H^{\frac 12}}\le |t|\max_{|s-t_0|\le |t|}\|u_s-u_{t_0}\|_{H^{\frac 12}},$$
which implies that
$$\lim_{t\to 0}\left\|\frac{U_{t_0+t}-U_{t_0}}{t}-u_{t_0}\right\|_{H^{\frac 12}}\le\lim_{t\to 0}\left(\max_{|s-t_0|\le |t|}\|u_s-u_{t_0}\|_{H^{\frac 12}}\right)=0,$$
that is, $U_t$ is differentiable at $t_0$, and (7.6) holds.\quad$\square$
\enddemo
\noindent {\bf Remark 7.5}\quad Here it is an appropriate place to relate  a result of Figalli [Fi]. In an attempt to study the regularity of the elements in $\WP(S^1)$, Figalli [Fi] investigated the smoothness  of the flows of the $H^{\frac 32}$ vector fields and showed that there exists some $H^{\frac 32}$  vector field whose flow  is neither Lipschitz nor $H^{\frac 32}$. Now our Theorem 7.3 says  that the the flow maps of the $H^{\frac 32}$ vector field in Figalli's example must also belong to $\WP(S^1)$, which in turn implies (Theorem 1.2) that there exists some quasisymmetric homeomorphism which is in $\WP(S^1)$ but is  neither  $H^{\frac 32}$ nor Lipschitz. Our proof of Theorem 1.2  relies neither  on Theorem 7.3 nor on Figalli's result. Moreover, it gives an explicit expression of  a quasisymmetric homeomorphism of the Weil-Petersson class being neither  $H^{\frac 32}$ nor Lipschitz.

\head 8. Proof of Theorem 1.4
\endhead

Recall that  the universal Teichm\"uller space has a quasisymmetric homeomorphism model, namely, $T=\QS(S^1)/\text{M\"ob}(S^1)$. Now $\Cal T=\QS(S^1)/\Rot(S^1)$ is a fiber space over $T$ and in fact is a model of the universal Teichm\"uller curve (see [Ber], [Te], [TT2]). Each point in $\Cal T$ can be considered as a quasisymmetric homeomorphism which keeps  1 fixed. There exists a one-to-one map $\Psi$ from $\Cal T$ onto $\hat T_b$ (another model of the universal Teichm\"uller curve)  which sends  $h$  to $\log f'$ under the normalized decomposition $h^{-1}=f^{-1}\circ g$. Via $\Psi$, $\Cal T$ is endowed with a standard complex Banach manifold structure such that $\Psi:\Cal T\to \hat T_b$ is a bi-holomorphic isomorphism (see [TT2] for more details).

Now we consider the Weil-Petersson class. Set $\Cal T_0=\WP(S^1)/\Rot(S^1)$. Then  $\Psi$
establishes a bijective map between $\Cal T_0$ and $\hat T_b\cap\Cal{AD}_0(\Delta)$. As stated in  Remark 5.1, a natural metric assigned to $\Cal T_0$ is the
following $H^{\frac 12}$ metric:
$$d(h_1, h_2)=\|\log |h'_2|-\log |h'_1|\|_{H^{\frac 12}},\quad
h_1,\,h_2\in \Cal T_0.\tag 8.1$$
Examining the last step in the proof of Theorem 1.1, we see that the metric is topologically equivalent to the following metric:
$$d'(h_1, h_2)=\|\log h'_2-\log h'_1\|_{H^{\frac 12}},\quad
h_1,\,h_2\in \Cal T_0.\tag 8.2$$

Then we have the following result.

\proclaim{Theorem 8.1} $\Psi: (\Cal T_0, d)\to\hat T_b\cap\Cal{AD}_0(\Delta)$ is a homeomorphism.
\endproclaim
\demo{Proof} We first recall the fact that $\|\log h'\|_{H^{\frac 12}}$ is small if and only if $\|\log (h^{-1})'\|_{H^{\frac 12}}$ is small. Examining the proof of Theorem 1.1, we find out that $\|\Psi(h)\|_{\Cal {AD}}$ is small if $\|\log h'\|_{H^{\frac 12}}$ is small. Thus, $\Psi$ is continuous at the base point $\id$.  Conversely, suppose $\log f'\in \hat T_b\cap\Cal{AD}_0(\Delta)$ has small norm. Let $h^{-1}=f^{-1}\circ g$ be the normalized conformal sewing mapping of $f$. We need to show that $\|\log h'\|_{H^{\frac 12}}$ is small, or equivalently, $\|\log (h^{-1})'\|_{H^{\frac 12}}$ is small.

Since $S_f=\Lambda(\log f')$ has small norm (2.5), by means of the well-known Ahlfors-Weil section (see [AW]),  $f$
can be extended to a quasiconformal mapping in the whole plane whose
complex dilatation $\mu$ has the form
$$\mu(z)=-\frac {1}{2}(|z|^2-1)^2S_f(\bar z^{-1})\bar
z^{-4},\, z\in\Delta^*.\tag 8.3$$ Thus, $\mu\in\Cal M(\Delta^*)$ with small norm $\|\mu\|_{WP}$. By means of Lemma 1.5 in [TT2], we have $f_{\mu}(\infty)=\infty$.

 We first consider the special case that $f=f_{\mu}|_{\Delta}$. Let $w_{\mu}$ be the unique quasiconformal mapping of $\Delta^*$ onto itself with Beltrami coefficient $\mu$ and keeping the points $1$ and $\infty$ fixed. Extending $w_{\mu}$ to the unit disk by symmetry, we obtain a quasiconformal mapping $w_{\mu}$ in the whole plane with $w_{\mu}(0)=0$. Then $g=f_{\mu}\circ w^{-1}_{\mu}|_{\Delta^*}$, and $h=w_{\mu}|_{S^1}$. Now lemma 2.5 in [TT2] implies that the Beltrami coefficient $\nu$ of $w^{-1}_{\mu}$ has small norm $\|\nu\|_{WP}$.  On the other hand, it is easy to see that  $h={g}^{-1}\circ  f$ is the  quasisymmetric conformal
 sewing mapping corresponding to $rj\circ g\circ j$, where $j(z)=\bar z^{-1}$ is  the standard reflection of the unit circle, and $r$ is a constant such that $ r(j\circ g\circ j)'(0)=1$.
 Now  $rj\circ g\circ j=rj\circ f_{\mu}\circ w^{-1}_{\mu}\circ j|_{\Delta}$ has the quasiconformal extension $rj\circ f_{\mu}\circ w^{-1}_{\mu}\circ j|_{\Delta^*}$ which keeps the point at infinity fixed, we conclude that  $\log(rj\circ g\circ j)'$ has small norm in $\Cal {AD}_0(\Delta)$ since the Beltrami coefficient $\nu$ of $w^{-1}_{\mu}$ has small norm $\|\nu\|_{WP}$. Thus, $\log g'$ has small norm in $\Cal{AD}(\Delta^*)$. It follows from (5.1) that $\|\log (h^{-1})'\|_{H^{\frac 12}}$ is small.

 In the general case, since $f$ and $f_{\mu}$ have the same complex dilatation $\mu$, we conclude by the normalized conditions $f(0)=f_{\mu}(0)=0$, $f'(0)=f'_{\mu}(0)=1$ and $f^{''}_{\mu}(0)=0$ that $f=\gamma_1\circ f_{\mu}$, where $\gamma_1(z)=\frac{z}{1-\lambda z}$ with $\lambda=f^{''}(0)/2$. Since $\log f'\in \hat T_b\cap\Cal{AD}_0(\Delta)$ has small norm, we conclude that $\lambda=f^{''}(0)/2$ is  small (see [Te]). To find the normalized conformal sewing map of $f_{\mu}$, we set $z_1=g^{-1}(-\frac{1}{\lambda})$, $z_2=-\frac{1-{z_1}}{{z_1}(1-\overline{z_1})}$, and
 $$\gamma_2(z)=\frac{1-\overline {z_2}}{1-z_2}\frac{z-z_2}{1-\overline{z_2}z}.$$
 A direct computation yields that $\gamma_2(1)=1$, $\gamma_2(\infty)=z_1$. Noting that $g(\infty)=\infty$, we conclude that $z_1$ tends to infinity and $z_2$ is small when $\lambda$ is small. Consider $\hat g=\gamma^{-1}_1\circ g\circ\gamma_2$. Then $\hat g$ is a conformal mapping from $\Delta^*$ onto $f_{\mu}(\Delta^*)$, and $\hat g(1)=f_{\mu}(1)$, $\hat g(\infty)=\infty$.
   Consequently, the normalized conformal sewing map of $f_{\mu}$ is
   $${\hat h}^{-1}=f^{-1}_{\mu}\circ\hat g=f^{-1}_{\mu}\circ \gamma^{-1}_1\circ g\circ\gamma_2=f^{-1}\circ g\circ\gamma_2=h^{-1}\circ\gamma_2,$$
   which implies that $h=\gamma_2\circ\hat h$. By what we have proved in the first (special) case, we conclude that $\|\log {\hat h}'\|_{H^{\frac 12}}$ is small when $\|\log f'\|_{\Cal {AD}}$ is small. On the other hand, when $\|\log f'\|_{\Cal {AD}}$ is small, $z_2$ is small, which implies that
   $$\|\log \gamma'_2\|^2_{\Cal {AD}}=\frac{1}{\pi}\iint_{\Delta}\frac{4|z_2|^2}{|1-\overline{z_2}z|^2}dxdy=4\log\frac{1}{1-|z_2|^2}$$
   is also small. Therefore, we conclude by $\log h'=\log\gamma'_2\circ\hat h+\log{\hat h}'$ that $\|\log {h}'\|_{H^{\frac 12}}$ is small when $\|\log f'\|_{\Cal {AD}}$ is small. This completes the proof  that
   $\Psi^{-1}$ is continuous at the base point $0$.

We now handle the general case by changing a general point to the base point. We only sketch the standard procedure by using the so-called allowable mappings (see [Ber], [Na], [TT2] for more details).
Let $h\in \Cal T_0$ be fixed. Consider the map $R_h$ defined by
$R_h(k)=k\circ h^{-1}$. Then $R_h$ is a bijective map from $\Cal T_0$
onto itself. Noting that
$$d'(R_h(k_1), R_h(h_2))=\|(\log k'_2-\log k'_1)\circ h^{-1}\|_{H^{\frac 12}},\tag 8.4$$
we conclude that $R_h$ is a quasi-isometric
map from $\Cal T_0$ onto itself under the $d'$-metric. Now let $w$ be a quasiconformal extension of $h$ to $\Delta^*$ such
that $w$ is  quasi-isometric under the Poincar\'e with  Beltrami
coefficient $\mu\in\Cal M(\Delta^*)$. The existence of such a quasiconformal extension is guaranteed by means of the well-known Douady-Earle [DE] extension of a quasisymmetric homeomorphism (see [Cu]). As stated in Proposition 2.2, $R_{w}$ induces a bi-holomorphic
isomorphism $w^*$ from $T_0$ onto itself with $w^*\circ\Phi=\Phi\circ R_w$. In fact, it is known that $R_{w}$ also induces a bi-holomorphic
isomorphism $\tilde {w^*}$ from $\hat T_b\cap\Cal{AD}_0(\Delta)$ onto itself which  is related to $R_h$ by $ \tilde {w^*}\circ\Psi=\Psi\circ R_h$.  By
using the allowable mappings $\tilde {w^*}$ and $R_h$, we conclude that both
$\Psi: (\Cal T_0, d)\to\hat T_b\cap\Cal{AD}_0(\Delta)$ and its inverse are continuous at a general point
$h$ (or $\Psi(h)$).\quad$\square$

\enddemo
\vskip 0.3 cm

\noindent {\bf Proof of Theorem 1.4:}\quad Let $\tilde\Lambda$ denote the natural projection from $\Cal T_0=\WP(S^1)/\Rot(S^1)$ onto $T_0=\WP(S^1)/\text{M\"ob}(S^1)$. The metric $d$ on $\Cal T_0$  descends down to a metric on $T_0$, still denoted by $d$, as follows:
$$d(h_1, h_2)=\inf\{d(\tilde h_1, \tilde h_2): \tilde\Lambda(\tilde h_1)=h_1, \tilde\Lambda(\tilde h_2)=h_2\},\quad
h_1,\,h_2\in T_0.\tag 8.5$$ By Theorem 8.1, $\Psi$ establishes a homeomorphism from $(\Cal T_0, d)$ onto $\hat T_b\cap\Cal{AD}_0(\Delta)$. On the other hand, Theorem 2.5 says  that $\Lambda: \hat T_b\cap\Cal{AD}_0(\Delta)\to\beta(T_0)$ is a holomorphic split submersion. This already implies that $\Psi: (\Cal T_0, d)\to\hat T_b\cap\Cal{AD}_0(\Delta)$ induces a homeomorphism from $\tilde \Psi: (T_0, d)$ onto $\beta(T_0)$, which implies that the metric $d$ and Weil-Petersson metric induce the same topology on  $T_0=\WP(S^1)/\text{M\"ob}(S^1)$.\quad$\square$

\head 9 Open problems
\endhead
It is known that $\Cal T_0=\WP(S^1)/\Rot(S^1)$ inherits a standard complex Hilbert manifold structure from $\Cal {AD}_0(\Delta)$ by the bijection  $\Psi: \Cal T_0\to\hat T_b\cap\Cal{AD}_0(\Delta)$ (see [TT2]). Meanwhile, $H^{\frac 12}_{\Bbb R}/\Bbb R$ provides $\Cal T_0$ with a real Hilbert manifold structure by the correspondence $h\mapsto\log|h'|$ (see Remark 5.1). Now Theorem 8.1 says that theses two Hilbert manifold structures induce the same topology on $\Cal T_0$.  It is not clear whether these two manifold structures are well compatible with each other. We believe this is the case and propose the following

\proclaim {Conjecture 9.1} Under the normalized decomposition $h^{-1}=f^{-1}\circ g$, both the bijective map $\log f'\mapsto\log|h'|$ and its inverse are real analytic. \endproclaim

In the recent paper [GR], Gay-Balmaz-Ratiu made the following

\proclaim {Conjecture 9.2 [GR]} Given a continuous vector field $u(t, \cdot)\in C^0([0, M], H^{\frac 32})$ with the normalized condition $(1.4)$,  the flow maps $h(t, \zeta)$ of the differential equation $$
\cases
\frac{dh}{dt}=u(t, h)\\
h(0,\zeta)=\zeta
\endcases
\tag 9.1$$
belong to the Weil-Petersson class, namely, $h(t, \cdot)\in \WP(S^1)$ for each fixed $t\in [0, M]$; Furthermore, the mapping $t\mapsto h(t, \cdot)$ from $[0, M]$ into $\WP(S^1)$ is continuously differentiable under the standard Hilbert manifold structure introduced by Takhtajan-Teo [TT2]. \endproclaim

The first assertion in Conjecture 9.2 is true by our Theorem 7.3. Furthermore, Theorem 7.3  implies that the mapping $t\mapsto \log|h'(t, \cdot)|$ from $[0, M]$ into $H^{\frac 12}_{\Bbb R}$ is continuously differentiable. It is clear that if Conjecture 9.1 were true, then Conjecture 9.2 would also be true.

Based on Lemma 7.1, it is natural to propose the following

\proclaim{Problem 9.3} Given a continuous vector field $u(t, \cdot)\in C^0([0, M], \Lambda^*)$ with the normalized condition $(1.4)$ and $(1.5)$, let $h(t, \cdot)$be  the flow maps of the differential equation $(9.1)$. Determine whether or not the flow $h(t, \cdot): [0, M]\rightarrow T$  is continuously differentiable.
\endproclaim

\head 10 Appendix: proof of Propositions 4.1 and 4.3
\endhead
In this section, we will prove Propositions 4.1 and 4.3 as we promised in Section 4. We restate them as follows.
\proclaim{Proposition 10.1} Let $h$ and $h_0$ be quasisymmetric homeomorphisms which keep the points $1$, $-1$ and $i$ fixed. Then for each fixed $u\in H^{\frac 12}$, $\|P_hu-P_{h_0}u\|_{H^{\frac 12}}\to 0$ when $\tau(h, h_0)\to 0$.
\endproclaim

\proclaim {Proposition 10.2} $P^+_h$ is a bounded isomorphism from $\Cal{AD}(\Delta)$ onto itself.  Moreover, it holds that
$$\|P^+_h\phi\|^2_{\Cal{AD}}=\|\phi\|^2_{\Cal{AD}}+\|P^-_h\phi\|^2_{\Cal{AD}},\,\phi\in\Cal{AD}(\Delta).\tag 10.1$$
\endproclaim

Here it is an appropriate place to point out that, though not stated in this form, Proposition 10.1 has appeared in the unpublished Master thesis [Li] of Q. Liu.
To prove Propositions 10.1 and 10.2, we need two related operators. Let $\Cal A^2(\Delta)$ denote the complex Hilbert space of all holomorphic
functions $\psi$ on the unit disk with norm
$$\|\psi\|_{\Cal A^2}=\left(\frac{1}{\pi}\iint_{\Delta}|\psi(\zeta)|^2d\xi d\eta\right)^{\frac
12}.\tag 10.2$$
Then $D\phi(z)=\phi'(z)$ determines an  isometric
isomorphism from $\Cal {AD}_0(\Delta)$ onto $\Cal A^2(\Delta)$.

For a quasisymmetric homeomorphism $h$, two
 kernel functions were introduced in the previous paper [HS] by Hu and
 the  author. They are
$$\phi_h(\zeta, z)=\frac{1}{2\pi
i}\int_{S^1}\frac{h(w)}{(1-\zeta w)^2(1-zh(w))}dw, \quad (\zeta,
z)\in\Delta\times\Delta,\tag 10.3$$
$$\psi_h(\zeta, z)=\frac{1}{2\pi
i}\int_{S^1}\frac{h(w)}{(\zeta- w)^2(1-zh(w))}dw, \quad (\zeta,
z)\in\Delta\times\Delta. \tag 10.4$$
The two kernels $\phi_h$ and $\psi_h$ induce two bounded operators on $\Cal A^2(\Delta)$ as follows:
$$T^-_h\psi(\zeta)=\frac{1}{\pi}\iint_{\Delta}\phi_h(\zeta, \bar  z)\psi(z)dxdy,\quad \psi\in\Cal A^2(\Delta),\, \zeta\in\Delta,\tag 10.5$$
and
$$T^+_h\psi(\zeta)=\frac{1}{\pi}\iint_{\Delta}\psi_h(\zeta, \bar z)\psi(z)dxdy,\quad \psi\in\Cal A^2(\Delta),\, \zeta\in\Delta.\tag 10.6$$
Then, Theorem 3.1 in [HS] says that on $\Cal {AD}(\Delta)$,
$$D\circ P^-_h=T^-_h\circ D,\quad D\circ P^+_h=T^+_h\circ D,\tag 10.7$$
while Lemma 2.3 in [SW] says that
$$\|T^+_h\psi\|^2_{\Cal A^2}=\|\psi\|^2_{\Cal A^2}+\|T^-_h\psi\|^2_{\Cal A^2},\quad \psi\in\Cal A^2(\Delta).\tag 10.8
$$

We first prove
\proclaim{Lemma 10.3} Let $h$ be a quasisymmetric homeomorphism which keep the points $1$, $-1$ and $i$ fixed. Then

$(1)$ $\|P^-_h\|\to 0$ when $\tau(0, h)\to 0$.

$(2)$ For each fixed $\phi\in\Cal {AD}$, $\|P^+_h\phi-\phi\|_{\Cal {AD}}\to 0$ as $\tau(0, h)\to 0$.

\endproclaim
\demo{Proof} By the definition (2.1) of the Teichm\"uller metric, there exists a so-called extremal quasiconformal extension $f$ of $h$ so that its complex dilatation $\mu$ satisfies
$$\|\mu\|_{\infty}=\frac{e^{2\tau(0, h)}-1}{e^{2\tau(0, h)}+1}.$$
Thus, as $\tau(0, h)\to 0$, $\|\mu\|_{\infty}\to 0$. Since  $h$  keeps the points $1$, $-1$ and $i$ fixed, we conclude by Strebel's approximation theorem (see [St]) that $\partial f(z)\to 1$ for a.e. $z\in\Delta$,  and $f(z)\to z$ locally uniformly in $\Delta$.

(1) Proposition 3.1 in [HS] says that
$$\|T^-_h\|\le\frac{\|\mu\|_{\infty}}{\sqrt{1-\|\mu\|^2_{\infty}}},\tag 10.9$$
which implies that, when $\tau(0, h)\to 0$, $\|T^-_h\|\to 0$ and consequently that $\|P^-_h\|\to 0$  by (10.7).

(2) By (10.7) we need to show that for each fixed $\psi\in\Cal A^2$, $\|T^+_h\psi-\psi\|_{\Cal A^2}\to 0$ as $\tau(0, h)\to 0$. Clearly,
$$\|T^+_h\psi-\psi\|^2_{\Cal A^2}=\|T^+_h\psi\|^2_{\Cal A^2}+\|\psi\|^2_{\Cal A^2}-\frac{2}{\pi}\Re \iint_{\Delta}T^+_h\psi(\zeta)\overline{\psi(\zeta)}d\xi d\eta.\tag 10.10$$
Proposition 3.2 in [HS] says that
$$
T^+_h\psi(\zeta)=\frac{1}{\pi}\iint_{\Delta}\frac{{\partial}f(w)\psi(f(w))}{(1-\zeta
\bar w)^2}dudv,\tag 10.11$$
and
$$\|T^+_h\|\le\frac{1}{\sqrt{1-\|\mu\|^2_{\infty}}}.\tag 10.12$$
Then,
$$\align
\iint_{\Delta}T^+_h\psi(\zeta)\overline{\psi(\zeta)}d\xi d\eta&=\frac{1}{\pi}\iint_{\Delta}\left(\iint_{\Delta}\frac{{\partial}f(w)\psi(f(w))}{(1-\zeta
\bar w)^2}dudv\right)\overline{\psi(\zeta)}d\xi d\eta\\
&=\frac{1}{\pi}\iint_{\Delta}\left(\iint_{\Delta}\frac{\overline{\psi(\zeta)}}{(1-\zeta
\bar w)^2}d\xi d\eta\right){\partial}f(w)\psi(f(w))dudv\\
&=\iint_{\Delta}{\partial}f(w)\psi(f(w))\overline{\psi(w)}dudv.\tag 10.13
\endalign
$$
By (10.10-10.13), we only need to show that
$$\lim_{\|\mu\|_{\infty}\to 0}\iint_{\Delta}\left(|\psi|^2-\Re({\partial}f\psi(f)\overline{\psi})\right)=0.\tag 10.14$$
Noting that
$$0\le||\psi|^2-{\partial}f\psi(f)\overline{\psi}|+|\psi|^2-|{\partial}f\psi(f)\overline{\psi}|\le 2|\psi|^2,$$
we conclude by Lesbegue's dominated convergence theorem that
$$
\align
&\lim_{\|\mu\|_{\infty}\to 0}\iint_{\Delta}(||\psi|^2-{\partial}f\psi(f)\overline{\psi}|+|\psi|^2-|{\partial}f\psi(f)\overline{\psi}|)\\
&=\iint_{\Delta}\lim_{\|\mu\|_{\infty}\to 0}(||\psi|^2-{\partial}f\psi(f)\overline{\psi}|+|\psi|^2-|{\partial}f\psi(f)\overline{\psi}|)=0.
\endalign$$
On the other hand,
$$
\align
\left(\iint_{\Delta}|{\partial}f\psi(f)\overline{\psi}|\right)^2&\le\iint_{\Delta}|\psi|^2\iint_{\Delta}|\psi(f)|^2|\partial f|^2\\
&=\iint_{\Delta}|\psi|^2\iint_{\Delta}\frac{|\psi|^2}{1-|\mu(f^{-1})|^2}\\
&\le\frac{1}{1-\|\mu\|^2_{\infty}}\left(\iint_{\Delta}|\psi|^2\right)^2.\endalign
$$
Combining these two inequalities together we obtain
$$\lim_{\|\mu\|_{\infty}\to 0}\iint_{\Delta}||\psi|^2-{\partial}f\psi(f)\overline{\psi}|=0.\tag 10.15$$
Now (10.14) follows from (10.15) by noting
$$||\psi|^2-\Re({\partial}f\psi(f)\overline{\psi})|\le ||\psi|^2-{\partial}f\psi(f)\overline{\psi}|.$$\quad$\square$

\enddemo

\vskip 0.2 cm

\noindent {\bf Proof of Proposition 10.1}\quad Let $u=\phi+\overline\psi$ be given. Noting that
$$\align
P_hu(z)-u(z)&=P_h\phi(z)+\overline{P_h\psi(z)}-\phi(z)-\overline{\psi(z)}\\
&=P^+_h\phi(z)-\phi(z)+\overline{P^+_h\psi(z)}-\overline{\psi(z)}+P^-_h\phi(\bar z)+\overline{P^-_h\psi(\bar z)},\endalign$$
we conclude by Lemma 10.3 that $\|P_hu-u\|_{H^{\frac 12}}\to 0$ when $\tau(0, h)\to 0$. Consequently,
$$\|P_hu-P_{h_0}u\|_{H^{\frac 12}}=\|P_{h_0}(P_{h\circ h^{-1}_0}u-u)\|_{H^{\frac 12}}\le \|P_{h_0}\|\|P_{h\circ h^{-1}_0}u-u\|_{H^{\frac 12}}\to 0$$
when $\tau(0, h\circ h^{-1}_0)=\tau(h, h_0)\to 0$.\quad$\square$

\vskip 0.2 cm

To prove Proposition 10.2, we also need the so-called Grunsky operator.
Consider the normalized decomposition $h=f^{-1}\circ g$ as before. Set
$$U(f,\zeta,z)=\frac{f'(\zeta)f'(z)}{[f(\zeta)-f(z)]^2}-\frac 1{(\zeta-z)^2},\quad (\zeta, z)\in \Delta\times \Delta.\tag 10.16$$
Then $S_f(z)=-6U(f, z, z)$ is the Schwarzian derivative of $f$. $f$
determines the so-called Grunsky operator on $\Cal A^2(\Delta)$, defined as
$$G_f\psi(\zeta)=\frac{1}{\pi}\iint_{\Delta}U(f,
\zeta, \bar z){\psi(z)}dxdy.\tag 10.17$$ It is known that $G_f$ is a
bounded operator from $\Cal A^2(\Delta)$ into itself with  $\|G_f\|<1$ (see [Po1], [Sh], [TT2]). The following relation was proved  by the author and Wei [SW]:
$$T_h^+\circ G_f=J\circ T_h^-\circ J,\tag 10.18$$
where $J$ is the operator defined by $J\phi(z)=\overline{\phi(\bar z)}$ so that $J^2=$id, $J\circ D=D\circ J$.

\vskip 0.2 cm

\noindent {\bf Proof of Proposition 10.2}\quad (10.1) follows directly from (10.7) and (10.8). Now let $\psi\in \Cal {AD}(\Delta)$ be given. Choose  $\omega\in\Cal{AD}_0(\Delta)$ so that $D\omega=-G_fJD\psi$. By (10.18) it holds that
$$JT^-_hD\psi+T^+_hD\omega=JT^-_hD\psi-T^+_hG_fJD\psi=0.$$
By (10.7) we obtain
$$D(P^+_h\omega+JP^-_h\psi)=T^+_hD\omega+JT^-_hD\psi=0.$$
Then,
$$PP_h(\psi+\overline\omega)=P^+_h\psi+\overline{JP^-_h\psi}+\overline{P^+_h\omega}+JP^-_h\omega=P^+_h\psi+JP^-_h\omega+\overline{P^+_h\omega(0)}.$$
Set $\phi=P^+_h\psi+JP^-_h\omega+\overline{P^+_h\omega(0)}$. Then $\phi\in\Cal {AD}(\Delta)$, and $P_{h^{-1}}\phi=\psi+\overline\omega$. Consequently, $P^+_{h^{-1}}\phi=\psi$, and $P^+_{h^{-1}}$ is surjective. Replacing $h^{-1}$ with $h$, we conclude that $P^+_h$ is surjective. \quad$\square$

 \vskip 0.5 cm

\vskip 0.2 cm\noindent
{\bf Acknowledgements} \quad
The author would like to thank the  referee for a very careful reading of the manuscript and for several corrections.

\enddocument